\magnification=\magstep1
\overfullrule=0pt
\parskip=6pt
\baselineskip=15pt
\headline={\ifnum\pageno>1 \hss \number\pageno\ \hss \else\hfill \fi}
\pageno=1
\nopagenumbers
\hbadness=1000000
\vbadness=1000000

\input epsf

\centerline{\bf EXPLICIT CALCULATIONS OF TENSOR PRODUCT COEFFICIENTS FOR $E_7$}

\vskip 15mm

\centerline{\bf Meltem Gungormez} \centerline{Dept. Physics, Fac.
Science, Istanbul Tech. Univ.} \centerline{34469, Maslak, Istanbul,
Turkey } \centerline{e-mail: gungorm@itu.edu.tr}

\vskip 5mm

\centerline{\bf Hasan R. Karadayi} \centerline{Dept. Physics, Fac.
Science, Istanbul Tech. Univ.} \centerline{34469, Maslak, Istanbul,
Turkey } \centerline{e-mail: karadayi@itu.edu.tr}

\vskip 25mm

\centerline{\bf{Abstract}}

\noindent We propose a new method to calculate coupling coefficients of
$E_7$ tensor products. Our method is based on explicit use of $E_7$ characters
in the definition of a tensor product.

\noindent When applying Weyl character formula for $E_7$ Lie algebra, one needs to make
sums over 2903040 elements of $E_7$ Weyl group. To implement such enormous sums,
we show we have a way which makes their calculations possible. This will be accomplished by decomposing an $E_7$ character into 72 participating $A_7$ characters.

\eject

\vskip 3mm
\noindent {\bf{I.\ INTRODUCTION }}
\vskip 3mm

Let $G_7=E_7, A_7$ and $\Lambda$, $\Lambda^\prime$ be two dominant weights of $G_7$ where  $R(\Lambda)$ and $R(\Lambda^\prime)$ are corresponding irreducible representations.
For general terms, we follow the book of Humphreys {\bf [1]} as ever.

Tensor product of these two irreducible representations is defined by,

$$R(\Lambda) \otimes R(\Lambda^\prime) = R(\Lambda+\Lambda^\prime) +
\sum_{\lambda \in S(\lambda+ \lambda^\prime)} t(\lambda<\Lambda+\Lambda^\prime) \ \ \
R(\lambda) \eqno(I.1) $$
where $S(\lambda+ \lambda^\prime)$ is the set of $\Lambda+\Lambda^\prime$ subdominants and
$t(\lambda<\Lambda+\Lambda^\prime)$' s are tensor coupling coefficients. Though Steinberg formula is the best known way {\bf [2]}, a natural way to calculate tensor coupling coefficients is also to solve the equation

$$Ch(\Lambda) \otimes Ch(\Lambda^\prime) = Ch(\Lambda+\Lambda^\prime) +
\sum_{\lambda \in S(\lambda+ \lambda^\prime)} t(\lambda<\Lambda+\Lambda^\prime) \ \ \
Ch(\lambda) \eqno(I.2) $$

\noindent for tensor coupling coefficients.
$Ch(\lambda)$ here is the character of an irreducible representation R($\lambda$) which corresponds to a dominant weight $\lambda$ and it is defined by the famous
Weyl Character formula:
$$ Ch(\lambda^+) = {A(\lambda^{++}) \over A(\rho_{_{G_7}})} \  \eqno(I.3)  $$
where for a weight $\mu$ in general
$$ A(\mu) \equiv \sum_{\sigma \in W(G_7)} \ \epsilon(\sigma) \ e^{\sigma(\mu)}
  \eqno(I.4) $$
$W(G_7)$ is the Weyl Group of $G_7$ and each and every element $\sigma$ is the so-called Weyl reflection while $\epsilon(\sigma)$ denotes its sign and $e^{\sigma(\lambda^{++}})$'s here are
known as formal exponentials.Thoroughout this work, we assume $\lambda^{++}$ denotes a strictly dominant weight defined for a dominant $\lambda^+$ by
$$ \lambda^{++} \equiv \rho_{_{G_7}} + \lambda^+  \eqno(I.5) $$
where $\rho_{_{G_7}}$ is the Weyl vector of $G_7$.

The crucial fact here is that
$$ \|W(E_7)\| = 2903040  \eqno(I.6) $$
where $ \| {\cal S} \| $ denotes order of set $ {\cal S}$. It is easy to see then to implement the sum in (I.4) would not be realizable explicitly. We, instead, propose 72 specifically chosen Weyl reflections which give us $A_7$ dominant weights participating within
the same $E_7$ Weyl orbit $W(\Lambda^+)$ for any $E_7$ dominant weight $\Lambda^+$.
As it is shown in the next section, this makes the evaluation of (I.4) realizable for $E_7$ but in terms of 72 $A_7$ characters and hence easily implementable.

\vskip 3mm
\noindent {\bf{II.\ $A_7$ DECOMPOSITION OF $E_7$ LIE ALGEBRA  }}
\vskip 3mm

For $i=1,2 \dots ,7$, let $\lambda_i$'s and $\alpha_i$'s be
respectively the fundamental dominant weigths and simple roots of
$A_7$ Lie algebra with the following Dynkin diagram
\midinsert \epsfxsize=5cm \centerline{\epsfbox{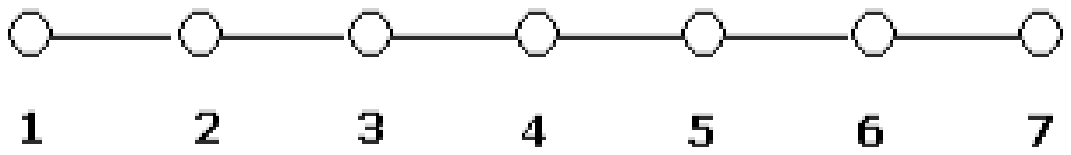}}
\endinsert
\noindent where $\rho_{_{A_7}} = \lambda_1 + \dots + \lambda_7 $  is $A_7$ Weyl vector and
$\Lambda_i$'s be fundamental dominant weights of $E_7$ Lie
algebra\ in according with the following Dynkin diagram,
\midinsert \epsfxsize=5cm \centerline{\epsfbox{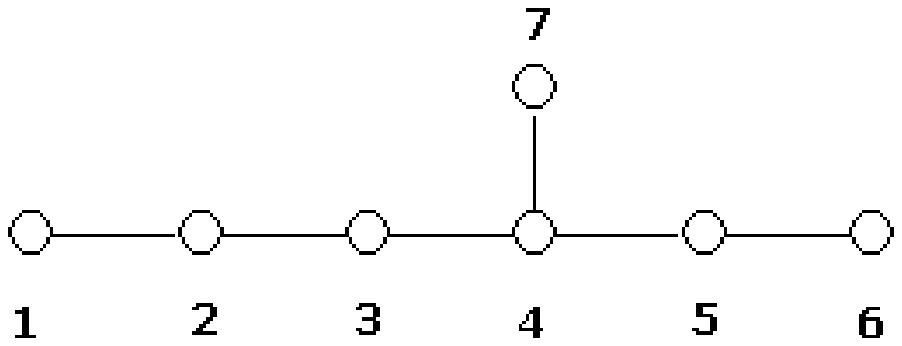}}
\endinsert
\noindent where $\rho_{_{E_7}} = \Lambda_1 + \dots + \Lambda_7 $ is $E_7$ Weyl vector.
\noindent We suggest following relations allows us to embed $A_7$ subalgebra into $E_7$ algebra:
$$ \eqalign{
\Lambda_1 &= \lambda_2  \cr
\Lambda_2 &= \lambda_1 + \lambda_3  \cr
\Lambda_3 &= 2 \ \lambda_3  \cr
\Lambda_4 &= 2 \ \lambda_3 + \lambda_6  \cr
\Lambda_5 &= \lambda_3 + \lambda_5  \cr
\Lambda_6 &= \lambda_4  \cr
\Lambda_7 &= \lambda_3 + \lambda_7  \cr }   \eqno(II.1) $$

This essentially means that
$$ {\| W(E_7)\| \over \|W(A_7)\| } = 72  \eqno(II.2) $$
\noindent which tells us that there are at most 72 $A_7$ dominant weights inside a
Weyl orbit $W(\Lambda^+)$. Note here that it is exactly 72 when $\Lambda^+$ is a strictly
dominant weight.
From the now on, $W(\mu)$ will always denotes the Weyl orbit of a weight $\mu$.

As the main point of view of this work, we present in appendix, 72 Weyl reflections to give 72 $A_7$ dominant weights participating in the same $E_7$ Weyl orbit $W(\Lambda^+)$ when they are exerted on the dominant weight $\Lambda^+$.
To this end, the Weyl reflections with respect to simple roots $\alpha_i$ will be called simple reflections
$\sigma_i$. We extend multiple products of simple reflections trivially by
$$ \sigma_{i_1,i_2}(\lambda) \equiv \sigma_{i_1}(\sigma_{i_2}(\lambda)) \ . $$
For $s=1, \dots 72$, $\Sigma(s)$'s are 72 Weyl reflections mentioned above.
As will also be seen by their definitions that,
$$\eqalign{
1) \ \ \epsilon(\sigma(s)) &=+1 \ \ \ s=1,2, \dots , 36 \cr
2) \ \ \epsilon(\sigma(s)) &=-1 \ \ \ s=37,38, \dots, 72 } $$
\vskip4mm
\noindent {\bf{III.\ CALCULATING TENSOR COUPLING COEFFICIENTS }}

\par Let us proceed in the instructive example

$$ R(\Lambda_3) \otimes R(\Lambda_4) = R(\Lambda_3+\Lambda_4) + \sum_{j=1}^{39} m(j) \
R(\theta_j) \eqno(III.1)  $$
of (I.1). One can see that there are 39 sub-dominant weights $ \theta_j $ of
$\Lambda_3+\Lambda_4 $ :
$$ \eqalign{
\theta_1 &= \Lambda_3 + \Lambda_4  \cr
\theta_2 &= \Lambda_1  \cr
\theta_3 &= \Lambda_7  \cr	
\theta_4 &= 3 \ \Lambda_1  \cr
\theta_5 &= \Lambda_1 + 2 \ \Lambda_2  \cr
\theta_6 &= 2 \ \Lambda_1 + \Lambda_3  \cr
\theta_7 &= \Lambda_1 + 2*\Lambda_5  \cr
\theta_8 &= \Lambda_1 + \Lambda_4 + \Lambda_6  \cr
\theta_9 &= 2 \ \Lambda_2 + \Lambda_7  \cr
\theta_{10} &= \Lambda_2 + \Lambda_3 + \Lambda_6  \cr
\theta_{11} &= \Lambda_2 + \Lambda_5 + \Lambda_7  \cr
} $$

$$ \eqalign{
\theta_{12} &= \Lambda_3 + 2 \ \Lambda_6  \cr
\theta_{13} &= 3 \ \Lambda_7  \cr	
\theta_{14} &= 3 \ \Lambda_1 + \Lambda_6  \cr
\theta_{15} &= \Lambda_1 + \Lambda_2 + 2 \ \Lambda_6 \cr
\theta_{16} &= \Lambda_1 + 3 \ \Lambda_6 \cr
\theta_{17} &= \Lambda_1 + \Lambda_6 + 2 \ \Lambda_7 \cr
\theta_{18} &= \Lambda_1 + \Lambda_2 + \Lambda_5  \cr
\theta_{19} &= \Lambda_1 + \Lambda_3 + \Lambda_7  \cr
\theta_{20} &= 2 \ \Lambda_1 + \Lambda_6 + \Lambda_7  \cr
\theta_{21} &= \Lambda_4 + \Lambda_7  \cr	
\theta_{22} &= \Lambda_5 + \Lambda_6 + \Lambda_7  \cr
\theta_{23} &= \Lambda_1 + \Lambda_2  \cr
\theta_{24} &= \Lambda_1 + \Lambda_6  \cr
\theta_{25} &= 2 \ \Lambda_1 + \Lambda_7  \cr	
\theta_{26} &= \Lambda_1 + 2 \ \Lambda_7  \cr	
\theta_{27} &= \Lambda_2 + \Lambda_3  \cr
\theta_{28} &= \Lambda_3  \cr
\theta_{29} &= \Lambda_3 + \Lambda_5  \cr
\theta_{30} &= 2 \ \Lambda_6 + \Lambda_7  \cr
\theta_{31} &= \Lambda_1 + \Lambda_5 + \Lambda_6  \cr
\theta_{32} &= \Lambda_1 + 2 \ \Lambda_6  \cr
\theta_{33} &= \Lambda_6 + \Lambda_7  \cr
\theta_{34} &= \Lambda_1 + \Lambda_4  \cr
\theta_{35} &= \Lambda_1 + \Lambda_5  \cr
\theta_{36} &= \Lambda_1 + \Lambda_2 + \Lambda_6  \cr
\theta_{37} &= \Lambda_2 + \Lambda_6 + \Lambda_7  \cr
\theta_{38} &= \Lambda_2 + \Lambda_7  \cr	
\theta_{39} &= \Lambda_3 + \Lambda_6  \cr
\theta_{40} &= \Lambda_5 + \Lambda_7	}$$

\noindent To this end, we should care about specialization of formal exponentials { \bf[3]}. Let us consider the so-called Fundamental Weights $\mu_I$  which are defined for $ (I=1, \dots ,8)$ as in the following { \bf[4]}:
$$ \alpha_i \equiv \mu_i - \mu_{i+1} \ \ , \ \ ( i=1, \dots 7.)  \eqno(III.2) $$
$\alpha_i$'s here are $A_7$ simple roots mentioned above and the best way to calculate $A_7$ and hence $E_7$ characters is to use the specialization in terms of parameters \ $ u_I \equiv e^{\mu_I}$ \ which are subjects of the condition \ $ \mu_1 + \mu_2 + \dots + \mu_8 = 0 $ \ or
\ $ u_1 u_2 \dots u_8 = 1 $.


To exemplify (I.3) for $E_7$, we would like to give detailed calculation of $Ch(\Lambda_3+\Lambda_4)$.
By applying 72 specifically chosen Weyl reflections on strictly dominant weight
$\rho_{_{E_7}}+\Lambda_3+\Lambda_4$, one can see we have the following decompositions:
$$ A(\rho_{_{E_7}}+\Lambda_3+\Lambda_4) = \sum_{k=1}^{36} Ch(\nu_k) - \sum_{k=37}^{72} Ch(\nu_k) \eqno(III.3) $$

where
$$\eqalign{
\nu_1 &=\lambda_1+\lambda_2+11 \ \lambda_3+\lambda_4+\lambda_5+2 \ \lambda_6+\lambda_7  \cr
\nu_2 &=2 \ \lambda_1+2 \ \lambda_2+8 \ \lambda_3+\lambda_4+\lambda_5+5 \ \lambda_6+\lambda_7  \cr
\nu_3 &=5 \ \lambda_1+\lambda_2 + 7 \ \lambda_3+\lambda_4+3 \ \lambda_5+2 \ \lambda_6+3 \ \lambda_7  \cr
\nu_4 &=4 \ \lambda_1+\lambda_2+6 \ \ \lambda_3+\lambda_4+3 \ \ \lambda_5+4 \ \ \lambda_6+
3 \ \lambda_7  \cr
\nu_5 &=2 \ \lambda_1+4 \ \lambda_2+4 \ \lambda_3+\lambda_4+5 \ \lambda_5+\lambda_6+
5 \ \lambda_7  \cr
\nu_6 &=5 \ \lambda_1+2 \ \lambda_2+5 \ \lambda_3+2 \ \lambda_4+2 \ \lambda_5+3 \ \lambda_6+
4 \ \lambda_7  \cr
\nu_7 &=5 \ \lambda_1+2 \ \lambda_2+5 \ \lambda_3+\lambda_4+4 \ \lambda_5+3 \ \lambda_6+
2 \ \lambda_7  \cr
\nu_8 &=7 \ \lambda_1+\lambda_2+6 \ \lambda_3+\lambda_4+2 \ \lambda_5+2 \ \lambda_6+
5 \ \lambda_7  \cr
\nu_9 &=7 \ \lambda_1+\lambda_2+5 \ \lambda_3+2 \ \lambda_4+3 \ \lambda_5+2 \ \lambda_6+
3 \ \lambda_7  \cr
\nu_{10} &=2 \ \lambda_1+2 \ \lambda_2+4 \ \lambda_3+\lambda_4+6 \ \lambda_5+\lambda_6+
6 \ \lambda_7  \cr
\nu_{11} &=3 \ \lambda_1+3 \ \lambda_2+3 \ \lambda_3+2 \ \lambda_4+4 \ \lambda_5+2 \ \lambda_6+
6 \ \lambda_7 \cr
\nu_{12} &=3 \ \lambda_1+3 \ \lambda_2+3 \ \lambda_3+\lambda_4+6 \ \lambda_5+2 \ \lambda_6+
4 \ \lambda_7  \cr
\nu_{13} &=6 \ \lambda_1+\lambda_2+5 \ \lambda_3+\lambda_4+2 \ \lambda_5+4 \ \lambda_6+
5 \ \lambda_7  \cr
\nu_{14} &=6 \ \lambda_1+\lambda_2+4 \ \lambda_3+2 \ \lambda_4+3 \ \lambda_5+4 \ \lambda_6+
3 \ \lambda_7  \cr
\nu_{15} &=\lambda_1+6 \ \lambda_2+\lambda_3+4 \ \lambda_4+2 \ \lambda_5+\lambda_6+
8 \ \lambda_7  \cr
\nu_{16} &=4 \ \lambda_1+4 \ \lambda_2+3 \ \lambda_3+\lambda_4+4 \ \lambda_5+\lambda_6+7 \ \lambda_7  \cr
\nu_{17} &=4 \ \lambda_1+4 \ \lambda_2+2 \ \lambda_3+2 \ \lambda_4+5 \ \lambda_5+\lambda_6+
5 \ \lambda_7  \cr
\nu_{18} &=\lambda_1+6 \ \lambda_2+\lambda_3+\lambda_4+8 \ \lambda_5+\lambda_6+2 \ \lambda_7 \cr
} $$
$$\eqalign{
\nu_{19} &=7 \ \lambda_1+2 \ \lambda_2+4 \ \lambda_3+\lambda_4+3 \ \lambda_5+3 \ \lambda_6+
4 \ \lambda_7  \cr
\nu_{20} &=8 \ \lambda_1+2 \ \lambda_2+2 \ \lambda_3+4 \ \lambda_4+\lambda_5+5 \ \lambda_6+\lambda_7  \cr
\nu_{21} &=10 \ \lambda_1+\lambda_2+3 \ \lambda_3+3 \ \lambda_4+\lambda_5+4 \ \lambda_6+
2 \ \lambda_7  \cr
\nu_{22} &=11 \ \lambda_1+\lambda_2+\lambda_3+6 \ \lambda_4+\lambda_5+2 \ \lambda_6+
\lambda_7  \cr
\nu_{23} &=\lambda_1+4 \ \lambda_2+\lambda_3+4 \ \lambda_4+3 \ \lambda_5+\lambda_6+
9 \ \lambda_7  \cr
\nu_{24} &=4 \ \lambda_1+2 \ \lambda_2+3 \ \lambda_3+\lambda_4+5 \ \lambda_5+\lambda_6+
8 \ \lambda_7  \cr
\nu_{25} &=4 \ \lambda_1+2 \ \lambda_2+2 \ \lambda_3+2 \ \lambda_4+6 \ \lambda_5+\lambda_6+
6 \ \lambda_7  \cr
\nu_{26} &=\lambda_1+4 \ \lambda_2+\lambda_3+\lambda_4+9 \ \lambda_5+\lambda_6+3 \ \lambda_7  \cr
\nu_{27} &=2 \ \lambda_1+5 \ \lambda_2+\lambda_3+3 \ \lambda_4+2 \ \lambda_5+2 \ \lambda_6+
9 \ \lambda_7  \cr
\nu_{28} &=5 \ \lambda_1+3 \ \lambda_2+2 \ \lambda_3+\lambda_4+5 \ \lambda_5+2 \ \lambda_6+
6 \ \lambda_7  \cr
\nu_{29} &=9 \ \lambda_1+\lambda_2+2 \ \lambda_3+3 \ \lambda_4+\lambda_5+6 \ \lambda_6+
2 \ \lambda_7  \cr
\nu_{30} &=\lambda_1+7 \ \lambda_2+\lambda_3+2 \ \lambda_4+2 \ \lambda_5+\lambda_6+
10 \ \lambda_7  \cr
\nu_{31} &=10 \ \lambda_1+2 \ \lambda_2+2 \ \lambda_3+2 \ \lambda_4+\lambda_5+
6 \ \lambda_6+\lambda_7  \cr
\nu_{32} &=13 \ \lambda_1+\lambda_2+\lambda_3+4 \ \lambda_4+\lambda_5+3 \ \lambda_6+
\lambda_7  \cr
\nu_{33} &=2 \ \lambda_1+2 \ \lambda_2+\lambda_3+5 \ \lambda_4+\lambda_5+\lambda_6+
12 \ \lambda_7  \cr
\nu_{34} &=\lambda_1+5 \ \lambda_2+\lambda_3+2 \ \lambda_4+3 \ \lambda_5+\lambda_6+
11 \ \lambda_7  \cr
\nu_{35} &=16 \ \lambda_1+\lambda_2+\lambda_3+2 \ \lambda_4+2 \ \lambda_5+\lambda_6+
\lambda_7  \cr
\nu_{36} &=\lambda_1+\lambda_2+3 \ \lambda_3+2 \ \lambda_4+\lambda_5+\lambda_6+
15 \ \lambda_7  \cr
\nu_{37} &=3 \ \lambda_1+\lambda_2+9 \ \lambda_3+\lambda_4+\lambda_5+4 \ \lambda_6+\lambda_7  \cr
\nu_{38} &=2 \ \lambda_1+\lambda_2+8 \ \lambda_3+\lambda_4+\lambda_5+6 \ \lambda_6+\lambda_7  \cr
\nu_{39} &=4 \ \lambda_1+2 \ \lambda_2+6 \ \lambda_3+\lambda_4+3 \ \lambda_5+3 \ \lambda_6+
3 \ \lambda_7  \cr
\nu_{40} &=6 \ \lambda_1+\lambda_2+6 \ \lambda_3+2 \ \lambda_4+2 \ \lambda_5+2 \ \lambda_6+
4 \ \lambda_7  \cr
\nu_{41} &=6 \ \lambda_1+\lambda_2+6 \lambda_3+\lambda_4+4 \lambda_5+2 \lambda_6+2 \lambda_7  \cr
\nu_{42} &=2 \ \lambda_1+3 \lambda_2+4 \lambda_3+\lambda_4+5 \lambda_5+2 \lambda_6+5 \lambda_7 \cr
\nu_{43} &=5 \ \lambda_1+\lambda_2+5 \lambda_3+2 \lambda_4+2 \lambda_5+4 \lambda_6+4 \lambda_7  \cr\
\nu_{44} &=5 \ \lambda_1+\lambda_2+5 \lambda_3+\lambda_4+4 \lambda_5+4 \lambda_6+2 \lambda_7  \cr
\nu_{45} &=3 \ \lambda_1+4 \lambda_2+3 \lambda_3+2 \lambda_4+4 \lambda_5+\lambda_6+6 \lambda_7  \cr
\nu_{46} &=3 \ \lambda_1+4 \lambda_2+3 \lambda_3+\lambda_4+6 \lambda_5+\lambda_6+4 \lambda_7  \cr
\nu_{47} &=6 \ \lambda_1+2 \lambda_2+5 \lambda_3+\lambda_4+2 \lambda_5+3 \lambda_6+5 \lambda_7  \cr
\nu_{48} &=6 \ \lambda_1+2 \lambda_2+4 \lambda_3+2 \lambda_4+3 \lambda_5+3 \lambda_6+3 \lambda_7  \cr
\nu_{49} &=8 \ \lambda_1+\lambda_2+5 \ \lambda_3+\lambda_4+3 \ \lambda_5+2 \ \lambda_6+
4 \ \lambda_7  \cr
\nu_{50} &=9 \ \lambda_1+\lambda_2+3 \ \lambda_3+4 \ \lambda_4+\lambda_5+4 \ \lambda_6+
\lambda_7 \cr
} $$
$$ \eqalign{
\nu_{51} &=3 \ \lambda_1+2 \ \lambda_2+3 \ \lambda_3+2 \ \lambda_4+5 \ \lambda_5+\lambda_6+
7 \ \lambda_7  \cr
\nu_{52} &=3 \ \lambda_1+2 \ \lambda_2+3 \ \lambda_3+\lambda_4+7 \ \lambda_5+\lambda_6+
5 \ \lambda_7  \cr
\nu_{53} &=\lambda_1+5 \ \lambda_2+\lambda_3+4 \ \lambda_4+2 \ \lambda_5+2 \ \lambda_6+
8 \ \lambda_7  \cr
\nu_{54} &=4 \ \lambda_1+3 \ \lambda_2+3 \ \lambda_3+\lambda_4+4 \ \lambda_5+2 \ \lambda_6+
7 \ \lambda_7 \cr
\nu_{55} &=4 \ \lambda_1+3 \ \lambda_2+2 \ \lambda_3+2 \ \lambda_4+5 \ \lambda_5+2 \ \lambda_6+
5 \ \lambda_7  \cr
\nu_{56} &=\lambda_1+5 \ \lambda_2+\lambda_3+\lambda_4+8 \ \lambda_5+2 \ \lambda_6+
2 \ \lambda_7  \cr
\nu_{57} &=7 \ \lambda_1+\lambda_2+4 \ \lambda_3+\lambda_4+3 \ \lambda_5+4 \ \lambda_6+
4 \ \lambda_7  \cr
\nu_{58} &=8 \ \lambda_1+\lambda_2+2 \ \lambda_3+4 \ \lambda_4+\lambda_5+6 \ \lambda_6+
\lambda_7  \cr
\nu_{59} &=2 \ \lambda_1+6 \ \lambda_2+\lambda_3+3 \ \lambda_4+2 \ \lambda_5+\lambda_6+
9 \ \lambda_7  \cr
\nu_{60} &=5 \ \lambda_1+4 \ \lambda_2+2 \ \lambda_3+\lambda_4+5 \ \lambda_5+\lambda_6+
6 \ \lambda_7  \cr
\nu_{61} &=9 \ \lambda_1+2 \ \lambda_2+2 \ \lambda_3+3 \ \lambda_4+\lambda_5+5 \ \lambda_6+
2 \ \lambda_7  \cr
\nu_{62} &=12 \ \lambda_1 + \lambda_2 + \lambda_3 + 5 \ \lambda_4 + \lambda_5 + 2 \ \lambda_6 +
2 \ \lambda_7  \cr
\nu_{63} &=11 \ \lambda_1+\lambda_2+3 \ \lambda_3+2 \ \lambda_4+\lambda_5+5 \ \lambda_6+
\lambda_7  \cr
\nu_{64} &=\lambda_1+2 \lambda_2+\lambda_3+6 \ \lambda_4+\lambda_5+\lambda_6+11 \ \lambda_7  \cr
\nu_{65} &=2 \ \lambda_1 + 4 \ \lambda_2 + \lambda_3 + 3 \ \lambda_4 + 3 \ \lambda_5 + \lambda_6 + 10 \ \lambda_7  \cr
\nu_{66} &=5 \ \lambda_1 + 2 \ \lambda_2 + 2 \ \lambda_3 + \lambda_4 + 6 \ \lambda_5 + \lambda_6 +7 \ \lambda_7  \cr
\nu_{67} &=\lambda_1 + 2 \ \lambda_2 + \lambda_3 + \lambda_4 + 11 \ \lambda_5 + \lambda_6 + \lambda_7  \cr
\nu_{68} &=\lambda_1 + 6 \ \lambda_2 + \lambda_3 + 2 \ \lambda_4 + 2 \ \lambda_5 + 2 \ \lambda_6 + 10 \ \lambda_7  \cr
\nu_{69} &=10 \ \lambda_1 + \lambda_2+2 \ \lambda_3+2 \lambda_4+\lambda_5+7 \lambda_6+\lambda_7  \cr
\nu_{70} &=15 \ \lambda_1 + \lambda_2 + \lambda_3+2 \ \lambda_4 + 3 \ \lambda_5 + \lambda_6 + \lambda_7  \cr
\nu_{71} &=\lambda_1+3 \ \lambda_2+\lambda_3+4 \ \lambda_4+\lambda_5+\lambda_6 +
13 \ \lambda_7  \cr
\nu_{72} &=\lambda_1+\lambda_2+2 \ \lambda_3+2 \ \lambda_4+\lambda_5+\lambda_6+16 \ \lambda_7   } $$
$A_7$ characters $Ch(\nu_k)$'s are defined by
$$ A(\rho_{_{A_7}}) \ Ch(\nu_k) = \sum_{\sigma \in W(A_7)} \ \epsilon(\sigma) \ e^{\sigma(\rho_{_{A_7}}+\nu_k)}  \eqno(III.4)  $$
\noindent Note here that $W(A_7)$ is the permutation group of 8 objects.

To display our result here, we use the following specialization of formal exponentials with only one free parameter x:
$$ \eqalign {
u_1 &= 1 \cr
u_2 &= 2 \cr
u_3 &= 3 \cr
u_4 &= 4 \cr
u_5 &= 5 \cr
u_6 &= 6 \cr
u_7 &= x \cr
u_8 &= 1/(720 \ x)  } \eqno(III.5) $$

In this specialization, one obtains the following one-parameter characters:
$$ \eqalign{
& A(\rho_{_{A_7}}) = -{1 \over 2^{20} \times 3^{11} \times 5^6 \ x^7 } \times \cr
& (-6 + x) \times (-5 + x) \times (-4 + x) \times (-3 + x) \times  (-2 + x) \times   \cr
& (-1 + x) \times (-1 + 720 \ x) \times (-1 + 1440 \ x) \times (-1 + 2160 \ x) \times  \cr
& (-1 + 2880 \ x) \times (-1 + 3600 \ x) \times (-1 + 4320 \ x) \times (-1 + 720 \ x^2) } $$
$$ \eqalign {
& A(\rho_{_{E_7}}) =  - { 7^3 \times 11 \times 13^2 \times 17^2 \times 23 \times 29 \times 47 \times 59^2 \times 71 \times 89 \times 179 \times 239 \times 359 \over
2^{40} \times 3^{19} \times 5^9 \times x^{10} }  \times \cr
& (-1 + 6 \ x) \times (-1 + 8 \ x) \times (-1 + 10 \ x) \times (-1 + 12 \ x)^2 \times (-1 + 15 \ x) \cr
& (-1 + 18 \ x) \times (-1 + 20 \ x) \times (-1 + 24 \ x)^2 \times (-1 + 30 \ x)^2 \times
(-1 + 36 \ x)  \cr
& (-1 + 40 \ x) \times (-1 + 48 \ x) \times (-1 + 60 \ x)^2 \times (-1 + 72 \ x) \times
(-1 + 90 \ x) \times (-1 +120 \ x)   } $$
$$ \eqalign{
& ch(\Lambda_4) = { 1 \over 2^{16} \times 3^8 \times 5^4 \ x^4 } \ (  \cr
& 2^2 \times 7 \times 29 \times 31  \times 113 \times 25849 \ +  \cr
& 3 \times 5^3 \times 7^2 \times 41 \times 13469 \times 25841 \ x +  \cr
& 2^2 \times 7^2 \times 227 \times 1997 \times 1004276389 \ x^2 +  \cr
& 2^2 \times 3^4 \times 5^2 \times 7^4 \times 17 \times 41 \times 13469 \times 45307 \ x^3 +  \cr
& 2^7 \times 3^2 \times 266944787316406807 \ x^4 +  \cr
& 2^6 \times 3^6 \times 5^3 \times 7^4 \times 17 \times 41 \times 13469 \times 45307 \ x^5 +  \cr
& 2^{10} \times 3^4 \times 5^2 \times 7^2 \times 227 \times 1997 \times 1004276389 \ x^6 +  \cr
& 2^{12} \times 3^7 \times 5^6 \times 7^2 \times 41 \times 13469 \times 25841 \ x^7 +  \cr
& 2^{18} \times 3^8 \times 5^4 \times 7 \times 29 \times 31 \times 113 \times 25849 \ x^8  )  } $$

$$ \eqalign {
& ch(\Lambda_3+\Lambda_4)={7 \over 2^{28} \times 3^{13} \times 5^6 \ x^7 } \ ( \cr
& 2^5 \ \ 7^2 \times 139 \times 40819 \times 22523219 + \cr
& 3 \times 296955329011336071883 \ \ x \ \ + \cr
& 7 \times 93629 \times 104327 \times 20612147800357 \ \ x^2 \ \ + \cr
& 2^2 \times 7 \times 89 \times 509 \times 407193532921684756441 \ \  x^3 \ \ + \cr
& 2^2 \times 7 \times 19\times 1447 \times 73091587 \times 1489316745532201 \ \  x^4 \ \ + \cr
& 2^4 \times 3 \times 11^2 \times 37 \times 1117 \times 18045889 \times 1661840436868789 \ \ x^5 \ \ + \cr
& 2^4 \times 3^2 \times 7 \times 11^2 \times 5237 \times 130069 \times 351401 \times 12440163841487 \ \  x^6 \ \ + \cr
& 2^9 \times 3^3 \times 11 \times 83 \times 757 \times 4830390973 \times 258355213888973 \ \ x^7 \ \ + \cr
& 2^8 \times 3^4 \times 5 \times 7 \times 11^2 \times 5237 \times 130069 \times 351401 \times 12440163841487 \ \ x^8 \ \ + \cr
& 2^{12} \times 3^5 \times 5^2 \times 11^2 \times 37 \times 1117 \times 18045889 \times 1661840436868789 \ \ x^9 \ \ + \cr
& 2^{14} \times 3^6 \times 5^3 \times 7 \times 19 \times 1447 \times 73091587 \times 1489316745532201 \ \ x^{10} \ \ + \cr
& 2^{18} \times 3^8 \times 5^4 \times 7 \times 89 \times 509 \times 407193532921684756441 \ \ x^{11} \ \ + \cr
& 2^{20} \times 3^{10} \times 5^5 \times 7  \times 93629 \times 104327 \times 20612147800357 \ \ x^{12} \ \ + \cr
& 2^{24} \times 3^{13} \times 5^6 \times 296955329011336071883 \ \ x^{13} \ \ + \cr
& 2^{33} \times 3^{14} \times 5^7 \times 7^2 \times 139 \times 40819 \times 22523219 \ \ x^{14}) \cr
& \ \ \ \ \cr
& ch(\theta_2)={1 \over 144 \ \ x} \ ( 3 \times 7 \times 17 + 5 \times 25841 \ x +
2^4 \times 3^3 \times 5 \times 7 \times 17 \ x^2 )  \cr
& \ \ \ \ \cr
&ch(\theta_3) = {1 \over 2^8 \times 3^4 \times 5 \ x^2 } \ (   \cr
& 5 \times 25841 + 3 \times 7 \times 19 \times 83 \times 2459 \ x +
2^2 \times 3^3 \times 7^2 \times 17 \times 45307 \ x^2 + \cr
& 2^4 \times 3^3 \times 5 \times 7 \times 19 \times 83 \times 2459 \ x^3 +
2^8 \times 3^4 \times 5^3 \times 25841 \ x^4 ) \cr
& \ \ \ \ \cr
&ch(\theta_4) = {1  \over 2^{12} \times 3^6 \times 5^2 \ x^3 } \ (   \cr
& 2 \times 3^2 \times 7 \times 251 \times 11287 \ +   \cr
& 2 \times 7 \times 11 \times 1634509733 \ x \ +  \cr
& 3^2 \times 7^2 \times 23 \times 73 \times 101 \times 281 \times 5783 \ x^2 \ + \cr
& 41 \times 953 \times 5987 \times 153547507 \ x^3 \ +  \cr
& 2^4 \times 3^4 \times 5 \times 7^2 \times 23 \times 73 \times 101 \times
281 \times 5783 \ x^4 \ + \cr
& 2^9 \times 3^4 \times 5^2 \times 7 \times 11 \times 1634509733 \  x^5 \ +  \cr
& 2^13 \times 3^8 \times 5^3 \times 7 \times 251 \times 11287 \ x^6 )   } $$

$$ \eqalign {
& ch(\theta_5) = {7 \over 2^{20} \times 3^{10} \times 5^4 \ x^5 } \ (   \cr
& 3^4 \times 389 \times 2621 \times 326611 \ +  \cr
& 44372305108670731 \ x \ +  \cr
& 2 \times 3 \times 7 \times 888351682544432651 \ x^2 \ +  \cr
& 2 \times 6835359680937443668841 \ x^3 \ +  \cr
& 2^3 \times 3 \times 17 \times 313 \times 2029 \times 4912124122975679 \ x^4 \ +  \cr
& 2^2 \times 7 \times 17 \times 87942983 \times 1353499222770007 \ x^5 \ +  \cr
& 2^7 \times 3^3 \times 5 \times 17 \times 313 \times 2029 \times
4912124122975679 \ x^6 \ +  \cr
& 2^9 \times 3^4 \times 5^2 \times 6835359680937443668841 \ x^7 \ +  \cr
& 2^{13} \times 3^7 \times 5^3 \times 7 \times 888351682544432651 \ x^8 \ +  \cr
& 2^{16} \times 3^8 \times 5^4 \times 44372305108670731 \ x^9 \ +  \cr
& 2^{20} \times 3^{14} \times 5^5 \times 389 \times 2621 \times 326611 \ x^{10})  } $$

$$ \eqalign {
& ch(\theta_6) = { 1 \over 2^{20} \times 3^{10} \times 5^4 \  x^5 } \ (  \cr
& 2 \times 3^5 \times 5 \times 7 \times 67 \times 118953227 \ +  \cr
& 2 \times 29 \times 191 \times 1259 \times 18041 \times 1205779 \ x \ +  \cr
& 3 \times 7 \times 419 \times 953 \times 23122934503549 \ x^2 +  \cr
& 7 \times 11 \times 13 \times 3881 \times 9624859 \times 1896757637 \ x^3 \ +  \cr
& 2^2 \times 3^2 \times 7 \times 111767 \times 14341763 \times 22885698919 \ x^4 \ +  \cr
& 2^2 \times 3^2 \times 5 \times 31 \times 114376828127 \times 589217525459 \ x^5 \ +  \cr
& 2^6 \times 3^4 \times 5 \times 7 \times 111767 \times 14341763 \times 22885698919 \ x^6 \ +  \cr
& 2^8 \times 3^4 \times 5^2 \times 7 \times 11 \times 13 \times 3881 \times 9624859 \times 1896757637 \ x^7 \ +  \cr
& 2^{12} \times 3^7 \times 5^3 \times 7 \times 419 \times 953 \times 23122934503549 \ x^8 \ +  \cr
& 2^{17} \times 3^8 \times 5^4 \times 29 \times 191 \times 1259 \times 18041 \times 1205779 \ x^9 \ +  \cr
& 2^{21} \times 3^{15} \times 5^6 \times 7 \times 67 \times 118953227 \ x^{10}   )   } $$
$$ \eqalign {
& ch(\theta_7) = { 7 \over 2^{26} \times 3^{14} \times 5^6 \ x^7 } \ (   \cr
& 2^2 \times 3^3 \times 71 \times 137 \times 897339719 \ +  \cr
& 2^3 \times 3^2 \times 9427181 \times 12037823647 \ x \ +  \cr
& 3 \times 5 \times 7^3 \times 6472827206445606859 \ x^2 \ +  \cr
& 7 \times 167 \times 1217 \times 3345931489 \times 3986554903 \ x^3 \ +  \cr
& 2^3 \times 3^2 \times 1913 \times 28712521177048302440611 \ x^4 \ +  \cr
& 2^3 \times 3^2 \times 5 \times 7 \times 31 \times 5504034331596437227418843 \ x^5 \ +  \cr
& 2^4 \times 3^5 \times 5^3 \times 157 \times 311 \times 6473 \times 6273427 \times
25424142019 \ x^6 \ +  \cr
& 2^{11} \times 3^5 \times 5 \times 11 \times 13 \times 193 \times 38358559 \times
321202718587369 \ x^7 \ +  \cr
& 2^8 \times 3^7 \times 5^4 \times 157 \times 311 \times 6473 \times 6273427 \times
25424142019 \ x^8 \ +  \cr
& 2^{11} \times 3^6 \times 5^3 \times 7 \times 31 \times
5504034331596437227418843 \ x^9 \ +  \cr
& 2^{15} \times 3^8 \times 5^3 \times 1913 \times 28712521177048302440611 \ x^{10} + \cr
& 2^{16} \times 3^8 \times 5^4 \times 7 \times 167 \times 1217 \times 3345931489 \times
3986554903 \ x^{11} +  \cr
& 2^{20} \times 3^{11} \times 5^6 \times 7^3 \times 6472827206445606859 \ x^{12} +  \cr
& 2^{27} \times 3^{14} \times 5^6 \times 9427181 \times 12037823647 \ x^{13} +  \cr
& 2^{30} \times 3^{17} \times 5^7 \times 71 \times 137 \times 897339719 \ x^{14} ) \cr
& \ \ \ \ \cr
& ch(\theta_8) = {7 \over 2^{26} \times 3^{12} \times 5^6 \ x^7 } \ (  \cr
& 70782069982080 \ +  \cr
& 2^5 \times 3  \times 7  \times 19  \times 101  \times 2437  \times 7883  \times 45307  \ x \ +  \cr
& 3^1  \times 839  \times 2591  \times 566506158976003  \ x^2 \ +  \cr
& 5^2 \times 127 \times 743526604818728544311 \ x^3 \ +  \cr
& 2^2  \times 5  \times 83  \times 263  \times 6791  \times 13901  \times 522761  \times 27245893 \ x^4 \ +  \cr
& 2^2 \times 1439  \times 6833 \times 57110377 \times 29440768222739  \ x^5 \ +  \cr
& 2^8 \times 3^2 \times 79 \times 2833 \times 5435185631 \times 1421253019763  \ x^6 \ +  \cr
& 2^6 \times 3^2 \times 5 \times 7 \times 37577587 \times 1652664511  \times 109347145543  \ x^7 \ +  \cr
& 2^{12} \times 3^4 \times 5 \times 79 \times 2833 \times 5435185631 \times 1421253019763  \ x^8 \ +  \cr
& 2^{10} \times 3^4 \times 5^2 \times 1439 \times 6833 \times 57110377 \times 29440768222739  \ x^9 \ +  \cr
& 2^{14}  \times 3^6 \times 5^4 \times 83 \times 263 \times 6791 \times 13901 \times 522761 \times 27245893  \  x^{10} \ +  \cr
& 2^{16} \times 3^8 \times 5^6 \times 127 \times 743526604818728544311 \ x^{11} \ +  \cr
& 2^{20} \times 3^{11} \times 5^5 \times 839 \times 2591 \times 566506158976003  \ x^{12} \ +  \cr
& 2^{29} \times 3^{13} \times 5^6 \times 7 \times 19 \times 101  \times 2437 \times 7883 \times 45307  \ x^{13} \ +  \cr
& 2^{35} \times 3^{15}  \times 5^8 \times 19 \times 101 \times 2437 \times 7883 \ x^{14} )
} $$

$$ \eqalign {
& ch(\theta_9) = { 1 \over 2^{23} \times 3^{12} \times 5^5 \ x^6 } \ (  \cr
& 5 \times 11 \times 1730263 \times 475374719 \ +  \cr
& 3 \times 5 \times 7 \times 11 \times 4933 \times 15316021986509 \ x \ +  \cr
& 5 \times 7 \times 13 \times 172831287782587727803 x^2 \ +  \cr
& 3 \times 7^2 \times 205820307753559239132761 x^3 \ +  \cr
& 2^2 \times 3^3 \times 41 \times 707606359 \times 1135814759505821 x^4 \ +  \cr
& 2 \times 3 \times 5 \times 7 \times 41 \times 1252496251 \times
20514266507074489 \ x^5 \ +  \cr
& 2^3  \times 3^3  \times 5  \times 67  \times 16903  \times 2819023  \times
2271886853005801  \  x^6 \ +  \cr
& 2^5  \times 3^3  \times 5^2  \times 7  \times 41  \times 1252496251  \times
20514266507074489  \ x^7 \ +  \cr
& 2^{10}  \times 3^7  \times 5^2  \times 41  \times 707606359  \times 1135814759505821  \  x^8 \ +  \cr
& 2^{12}  \times 3^7  \times 5^3  \times 7^2  \times 205820307753559239132761  \  x^9 \ +  \cr
& 2^{16}  \times 3^8  \times 5^5  \times 7  \times 13  \times 172831287782587727803  \ x^{10} \ +  \cr
& 2^{20}  \times  3^{11}  \times 5^6  \times 7  \times 11  \times 4933  \times 15316021986509  \ x^{11} \ +  \cr
& 2^{24}  \times  3^{12}  \times 5^7  \times 11  \times 1730263  \times 475374719  \ x^{12}  ) \cr
& \ \ \ \ \cr
& ch(\theta_{10}) = {1 \over 2^{26} \times 3^{12} \times 5^6 \ x^7 } \ (   \cr
& 2^2 \times 3 \times 5 \times 7^2 \times 63113 \times 7700299 \ +  \cr
& 3 \times 7^3 \times 45307 \times 63113 \times 7700299 \ x \ +  \cr
& 3 \times 5^2 \times 7 \times 17 \times 191 \times 30444055104860819 \ x^2 \ +  \cr
& 2^2 \times 13 \times 854020492656685314145399 \ x^3 \ +  \cr
& 2^2 \times 3 \times 7 \times 13 \times 137 \times 1619903 \times
40109150434199081 \ x^4 \ +  \cr
& 2^4 \times 3 \times 7 \times 13 \times 948542029 \times 239506095510065339 \ x^5 \ +  \cr
& 2^4 \times 3 \times 7 \times 11485393 \times 14627879671585167600299 \ x^6 \ +  \cr
& 2^6 \times 3^7 \times 5^2 \times 19 \times 61 \times 883 \times 8053 \times 44101 \times
244021 \times 6224627 \ x^7 \ +  \cr
& 2^8 \times 3^3 \times 5 \times 7 \times 11485393 \times
14627879671585167600299 \ x^8 \ + \  \cr
& 2^{12} \times 3^5 \times 5^2 \times 7 \times 13 \times 948542029 \times
239506095510065339 \ x^9 \ +  \cr
& 2^{14} \times 3^7 \times 5^3 \times 7 \times 13 \times 137 \times 1619903, 1 \times
40109150434199081 \ x^{10} \ +  \cr
& 2^{18} \times 3^8 \times 5^4 \times 13 \times 854020492656685314145399 \ x^{11}  \ +  \cr
& 2^{20} \times 3^{11} \times 5^7 \times 7 \times 17 \times 191 \times 30444055104860819 \ x^{12} \ +  \cr
& 2^{24} \times 3^{13} \times 5^6 \times 7^3 \times 45307 \times 63113 \times 7700299 \ x^{13} \ +  \cr
& 2^{30} \times 3^{15} \times 5^8 \times 7^2 \times 63113 \times 7700299 \ x^{14}   )     }  $$
$$ \eqalign {
& ch(\theta_{11}) = { 11 \over 2^{26} \times 3^{14} \times 5^6 \times x^7 } \ (   \cr
& 3^2 \times 7 \times 233 \times 2439503305003 +   \cr
& 251580659156174514557 \ x \ +   \cr
& 2 \times 3 \times 7 \times 197 \times 479 \times 83919415489238839 \ x^2 \ +   \cr
& 2 \times 3^5 \times 7 \times 257 \times 13163 \times 561199 \times 23985981091 \ x^3 \ + \cr  & 2^4 \times 3 \times 5^2 \times 7 \times 23 \times 31 \times 4297 \times 5779 \times
1802909 \times 106751927 \ x^4 \ +   \cr
& 2^3 \times 3^3 \times 207797403557 \times 58222225954444073 \ x^5 \ +   \cr
& 2^5 \times 3^3 \times 7 \times 22924263480995135201368307741 \ x^6 \ +   \cr
& 2^6 \times 3^5 \times 29 \times 10268423472246027440877336011 \ x^7 \ +   \cr
& 2^9 \times 3^5 \times 5 \times 7 \times 22924263480995135201368307741 \ x^8 \ +   \cr
& 2^{11} \times 3^7 \times 5^2 \times 207797403557 \times 58222225954444073 \ x^9 \ +   \cr
& 2^{16} \times 3^7 \times 5^5 \times 7 \times 23 \times 31 \times 4297 \times
5779 \times 1802909 \times 106751927 \ x^{10} \ +   \cr
& 2^{17} \times 3^{13} \times 5^4 \times 7 \times 257 \times 13163 \times 561199 \times
23985981091 \ x^{11} \ +   \cr
& 2^{21} \times 3^{11} \times 5^5 \times 7 \times 197 \times 479 \times
83919415489238839 \ x^{12} \ +   \cr
& 2^{24} \times 3^{12} \times 5^6 \times 251580659156174514557 \ x^{13} \ +   \cr
& 2^{28} \times 3^{16} \times 5^7 \times 7 \times 233 \times 2439503305003 \ x^{14}  ) \cr
& \ \ \ \ \cr
& ch(\theta_{12}) = {1 \over 2^{24} \times 3^{12} \times 5^6  \times x^7 } \ (   \cr
& 2^4 \times 3^3 \times 5^2 \times 7^2 \times 41 \times 13469 \ +    \cr
& 2^2 \times 3^3 \times 5 \times 7^3 \times 41 \times 13469 \times 45307  \ x \ +    \cr
& 3 \times 7^3 \times 41 \times 13469 \times 115259 \times 666353  \ x^2 \ +    \cr
& 5 \times 2999 \times 75703 \times 67759950443299 \  x^3 \ +    \cr
& 2^2 \times 3 \times 7 \times 41 \times 73^2 \times 367 \times
2736080070054181  \ x^4 \ +    \cr
& 2^2 \times 3^2 \times 7 \times 7685969351914603034407849  \ x^5 \ +    \cr
& 2^4 \times 3^3 \times 5^2 \times 7^2 \times 373 \times 571 \times 973787 \times
1054751382911  \ x^6 \ +    \cr
& 2^6 \times 3^4 \times 5 \times 37 \times 223 \times 311 \times 3767 \times
12319481 \times 1298488019  \ x^7 \ +    \cr
& 2^8 \times 3^5 \times 5^3 \times 7^2 \times 373 \times 571 \times 973787 \times
1054751382911 \ x^8 \  +    \cr
& 2^{10} \times 3^6 \times 5^2 \times 7 \times 7685969351914603034407849 \  x^9 \ +    \cr
& 2^{14} \times 3^7 \times 5^3 \times 7 \times 41 \times 73^2 \times 367 \times
2736080070054181 \ x^{10} \ +    \cr
& 2^{16} \times 3^8 \times 5^5 \times 2999 \times 75703 \times
67759950443299 \  x^{11} \ +    \cr
& 2^{20} \times 3^{11} \times 5^5 \times 7^3 \times 41 \times 13469 \times 115259 \times
666353  \ x^{12} \ +    \cr
& 2^{26} \times 3^{15} \times 5^7 \times 7^3 \times 41 \times 13469 \times
45307  \ x^{13} \ +    \cr
& 2^{32} \times 3^{17} \times 5^9 \times 7^2 \times 41 \times 13469 \  x^{14}   )  } $$

$$ \eqalign {
& ch(\theta_{13}) =  { 11 \over 2^{24} \times 3^{12} \times 5^5  \ x^6 } \ (  \cr
& 3238985337918907  +    \cr
& 3^4 \times 5 \times 7 \times 841610704258279  \ x  \ +    \cr
& 7^2 \times 13 \times 1482820942905285871  \ x^2  \ +    \cr
& 3^3 \times 7 \times 193 \times 467 \times 5690719 \times
2565345073  \ x^3 \ +    \cr
& 2^2 \times 3^2 \times 626862641769013044421771 \ x^4 \ +    \cr
& 2^4 \times 3^3 \times 5 \times 7 \times 79134994435751528887903  \ x^5 \ +    \cr
& 2^6 \times 3^4 \times 31 \times 149 \times 181 \times
9313661572349696183  \ x^6 \ +    \cr
& 2^8 \times 3^5 \times 5^2 \times 7 \times
79134994435751528887903  \ x^7 \ +    \cr
& 2^{10} \times 3^6 \times 5^2 \times 626862641769013044421771  \ x^8 \ +    \cr
& 2^{12}  \times 3^9 \times 5^3 \times 7 \times 193 \times 467 \times
5690719 \times 2565345073  \ x^9 \ +    \cr
& 2^{16} \times 3^8 \times 5^4 \times 7^2 \times 13 \times
1482820942905285871 \  x^{10}  \ +    \cr
& 2^{20} \times 3^{14} \times 5^6 \times 7 \times 841610704258279  \ x^{11} \ +  \cr
& 2^{24} \times 3^{12} \times 5^6 \times 3238985337918907  \ x^{12}   )   } $$

$$ \eqalign {
& ch(\theta_{14}) = { 1 \over 2^{18} \times 3^8 \times 5^4 \ x^5 } \ (  \cr
& 2^3 \times 3^2 \times 5 \times 7 \times 251 \times 11287 \ +    \cr
& 2 \times 3^2 \times 7^2 \times 251 \times 11287 \times 45307  \ x \ +    \cr
& 2 \times 7^2 \times 13 \times 47 \times 5077 \times 280370003  \ x^2 \ +    \cr
& 1305047 \times 32776442887619  \ x^3 \ +    \cr
& 7 \times 29 \times 67792850454653521231  \ x^4 \ +    \cr
& 2^2 \times 5 \times 59 \times 557 \times 739 \times 1838570856221483  \ x^5 \ +    \cr
& 2^4 \times 3^2 \times 5 \times 7 \times 29 \times 67792850454653521231  \ x^6 \ +    \cr
& 2^8 \times 3^4 \times 5^2 \times 1305047 \times 32776442887619  \ x^7 \ +    \cr
& 2^{13} \times 3^6 \times 5^3 \times 7^2 \times 13 \times 47 \times 5077 \times
280370003 \ x^8 \ +    \cr
& 2^{17} \times 3^{10} \times 5^4 \times 7^2 \times 251 \times
11287 \times 45307  \ x^9 \ +    \cr
& 2^{23} \times 3^{12} \times 5^6 \times 7 \times 251 \times 11287 \ x^{10}  )  } $$
$$ \eqalign {
& ch(\theta_{15}) =  { 7 \over 2^{24} \times 3^{12} \times 5^6 \times x^7 } \ (  \cr
& 2^{10} \times 3^3 \times 5^3 \times 51977 \ +    \cr
& 2^8 \times 3^3 \times 5^2 \times 7 \times 45307 \times 51977 \ x +  \ +    \cr
& 2^6 \times 3 \times 5 \times 7 \times 51977 \times 115259 \times 666353 \ x^2 \ +    \cr
& 3^2 \times 861317 \times 4060866137260903 \ x^3 \ +    \cr
& 3 \times 5 \times 41 \times 53 \times 16481 \times 24852500017071257 \ x^4 \ +    \cr
& 2^3 \times 223 \times 1035354314616125752647209 \ x^5 \ +    \cr
& 2^6 \times 3^2 \times 5 \times 495553698751 \times 90995390639929 \ x^6 \ +    \cr
& 2^7 \times 3^2 \times 5 \times 7 \times 163 \times 227 \times 99661 \times
33171175875735791 \ x^7 \ +    \cr
& 2^{10} \times 3^4 \times 5^2 \times 495553698751 \times 90995390639929 \ x^8 \ +    \cr
& 2^{11} \times 3^4 \times 5^2 \times 223 \times 1035354314616125752647209 \ x^9 \ +    \cr
& 2^{12} \times 3^7 \times 5^4 \times 41 \times 53 \times 16481 \times 24852500017071257 \ x^{10} \ +    \cr
& 2^{16} \times 3^{10} \times 5^4 \times 861317 \times 4060866137260903 \ x^{11} \ +    \cr
& 2^{26} \times 3^{11} \times 5^6 \times 7 \times 51977 \times 115259 \times
666353 \ x^{12} \ +    \cr
& 2^{32} \times 3^{15} \times 5^8 \times 7 \times 45307 \times 51977 \ x^{13} \ +    \cr
& 2^{38} \times 3^{17} \times 5^{10} \times 51977 \ x^{14}  )  \cr
& \ \ \ \ \cr
& ch(\theta_{16}) = { 1 \over 2^{22} \times 3^{10} \times 5^6 \times x^7 } \ (  \cr
& 2^6 \times 3^3 \times 5^3 \times 7 \times 17 \ +  \cr
& 2^4 \times 3^3 \times 5^2 \times 7^2 \times 17 \times 45307 \ x \ +  \cr
& 2^2 \times 3 \times 5 \times 7^2 \times 17 \times 115259 \times 666353 \ x^2 \ +  \cr
& 3 \times 7^3 \times 17 \times 37 \times 221303 \times 196292779 \ x^3 \ +  \cr
& 5 \times 7 \times 29 \times 147145853 \times 90252382399 \ x^4 \ +  \cr
& 2^2 \times 5 \times 86214774415429023276727 \ x^5 \ +  \cr
& 2^4 \times 3^4 \times 5^2 \times 7 \times 863 \times 4003 \times 33749 \times
4576169137 \ x^6 \ +  \cr
& 2^6 \times 3^3 \times 5^2 \times 23 \times 43 \times 7907 \times
14463096839915333 \ x^7 \ +  \cr
& 2^8 \times 3^6 \times 5^3 \times 7 \times 863 \times 4003 \times
33749 \times 4576169137 \ x^8 \ +  \cr
& 2^{10} \times 3^4 \times 5^3 \times 86214774415429023276727 \ x^9 \ +  \cr
& 2^{12} \times 3^6 \times 5^4 \times 7 \times 29 \times 147145853 \times
90252382399 \ x^{10} \ +  \cr
& 2^{16} \times 3^9 \times 5^4 \times 7^3 \times 17 \times 37 \times
221303 \times 196292779 \ x^{11} \ +  \cr
& 2^{22} \times 3^{11} \times 5^6 \times 7^2 \times 17 \times 115259 \times
666353 \ x^{12} \ +  \cr
& 2^{28} \times 3^{15} \times 5^8 \times 7^2 \times 17 \times
45307 \ x^{13} \ +  \cr
& 2^{34} \times 3^{17} \times 5^{10} \times 7 \times 17 \ x^{14}  )  } $$
$$ \eqalign {
& ch(\theta_{17}) =   { 7 \over 2^{26} \times 3^{12} \times 5^6 \times x^7 } \ (  \cr
& 2^2 \times 3^2 \times 5 \times 1920917020339 \ +  \cr
& 3^2 \times 7 \times 45307 \times 1920917020339 \ x \ +  \cr
& 2^2 \times 11 \times 1699 \times 20848837 \times 3497751463 \ x^2 \ +  \cr
& 2^3 \times 3 \times 13 \times 31 \times 47 \times 2517821 \times
2229533970737 \ x^3 \ +  \cr
& 5 \times 23 \times 239 \times 22442053352210346516283 \ x^4 \ +  \cr
& 2^3 \times 3 \times 11 \times 244186487395311739881589177 \ x^5 \ +  \cr
& 2^8 \times 3^3 \times 11^2 \times 3229 \times 1352409142907613715103 \ x^6 \ +  \cr
& 2^8 \times 3^2 \times 5 \times 11 \times 993682140710634056756709367 \ x^7 \ +  \cr
& 2^{12} \times 3^5 \times 5 \times 11^2 \times 3229 \times
1352409142907613715103 \ x^8 \ +  \cr
& 2^{11} \times 3^5 \times 5^2 \times 11 \times
244186487395311739881589177 \ x^9 \ +  \cr
& 2^{12} \times 3^6 \times 5^4 \times 23 \times 239 \times
22442053352210346516283 \ x^{10} \ +  \cr
& 2^{19} \times 3^9 \times 5^4 \times 13 \times 31 \times 47 \times
2517821 \times 2229533970737 \ x^{11} \ +  \cr
& 2^{22} \times 3^{10} \times 5^5 \times 11 \times 1699 \times 20848837 \times
3497751463 \ x^{12} \ +  \cr
& 2^{24} \times 3^{14} \times 5^6 \times 7 \times 45307 \times
1920917020339 \ x^{13} \ +  \cr
& 2^{30} \times 3^{16} \times 5^8 \times 1920917020339 \ x^{14}  ) \cr
& \ \ \ \ \cr
& ch(\theta_{18}) = { 1 \over 2^{22} \times 3^{12} \times 5^5 \times x^6 } \ (  \cr
& 2^2 \times 3^3 \times 5 \times 11 \times 36319 \times 4281313 \ +  \cr
& 2^4 \times 3 \times 5 \times 7 \times 11^2 \times 46549 \times 679298969  \ x \ +  \cr
& 3^2 \times 5^2 \times 7 \times 11 \times 41 \times 10712225499148973  \ x^2 \ +  \cr
& 3 \times 5^2 \times 7 \times 569 \times 476209471 \times 24269673079  \ x^3 \ +  \cr
& 2^7 \times 5^2 \times 7 \times 717343331 \times 35775113699503  \ x^4 \ +  \cr
& 2^8 \times 3^2 \times 7 \times 20327 \times 3472171 \times 37024507616807  \ x^5 \ +  \cr
& 2^4 \times 3^2 \times 3301 \times 3253233143634211446037993  \ x^6 \ +  \cr
& 2^{12} \times 3^4 \times 5 \times 7 \times 20327 \times 3472171 \times
37024507616807  \ x^7 \ +  \cr
& 2^{15} \times 3^4 \times 5^4 \times 7 \times 717343331 \times
35775113699503  \ x^8 \ +  \cr
& 2^{12} \times 3^7 \times 5^5 \times 7 \times 569 \times
476209471 \times 24269673079  \ x^9 \ +  \cr
& 2^{16} \times 3^{10} \times 5^6 \times 7 \times 11 \times 41 \times
10712225499148973  \ x^{10} \ +  \cr
& 2^{24} \times 3^{11} \times 5^6 \times 7 \times 11^2 \times
46549 \times 679298969  \ x^{11} \ +  \cr
& 2^{26} \times 3^{15} \times 5^7 \times 11 \times 36319 \times
4281313  \ x^{12}  ) } $$
$$ \eqalign {
& ch(\theta_{19}) = { 7 \over 2^{24} \times 3^{12} \times 5^5 \ x^6 } \ (  \cr
& 11 \times 19 \times 54132425994601 \ +  \cr
& 2^2 \times 3 \times 1873 \times 801641 \times 1539041989  \ x \ +  \cr
& 2^3 \times 5 \times 3391 \times 157884472828805717 \ x^2 \ +  \cr
& 3 \times 5 \times 41 \times 421 \times 26894683525927701469 \ x^3 \ +  \cr
& 2^3 \times 5296547 \times 40807969 \times 561290628359 \ x^4 \ +  \cr
& 2^2 \times 3^2 \times 5 \times 7 \times 18542831 \times 2641233933522645623 \ x^5 \ +  \cr
& 2^4 \times 3^4 \times 11699291 \times 71109193 \times 2011720082303 \ x^6 \ +  \cr
& 2^6 \times 3^4 \times 5^2 \times 7 \times 18542831 \times 2641233933522645623 \ x^7 \ +  \cr
& 2^{11} \times 3^4 \times 5^2 \times 5296547 \times 40807969 \times 561290628359 \ x^8  +  \cr
& 2^{12} \times 3^7 \times 5^4 \times 41 \times 421 \times 26894683525927701469 \ x^9 \ +  \cr
& 2^{19} \times 3^8 \times 5^5 \times 3391 \times 157884472828805717 \ x^{10} \ +  \cr
& 2^{22} \times 3^{11} \times 5^5 \times 1873 \times 801641 \times 1539041989 \ x^{11} \ +  \cr
& 2^{24} \times 3^{12} \times 5^6 \times 11 \times 19 \times 54132425994601 \ x^{12}  ) } $$

$$ \eqalign {
& ch(\theta_{20}) = { 7 \over 2^{22} \times 3^{10} \times 5^5 \times x^6 } \ (  \cr
& 2^3 \times 5 \times 11 \times 1634509733 +  \cr
& 2 \times 7 \times 11 \times 45307 \times 1634509733 \ x \ +  \cr
& 3^1 \times 5^2 \times 188673113 \times 883791547 \ x^2 \ +  \cr
& 2^5 \times 821 \times 5573 \times 41867209636027 \ x^3 \ +  \cr
& 3^4 \times 269 \times 27680957 \times 2625874917649 \ x^4 \ +  \cr
& 2^3 \times 3^2 \times 13 \times 29 \times 71 \times 74073603909992156759 \ x^5 \ +  \cr
& 2^5 \times 3^4 \times 5 \times 1489 \times 291506349305759397989 \ x^6 \ +  \cr
& 2^7 \times 3^4 \times 5 \times 13 \times 29 \times 71 \times
74073603909992156759  \ x^7 \ +  \cr
& 2^8 \times 3^8 \times 5^2 \times 269 \times 27680957 \times 2625874917649 \ x^8 \ +  \cr
& 2^{17} \times 3^6 \times 5^3 \times 821 \times 5573 \times 41867209636027  \ x^9 \ +  \cr
& 2^{16} \times 3^9 \times 5^6 \times 188673113 \times 883791547 \ x^{10} \ +  \cr
& 2^{21} \times 3^{10} \times 5^5 \times 7 \times 11 \times 45307 \times
1634509733 \ x^{11} \ +  \cr
& 2^{27} \times 3^{12} \times 5^7 \times 11 \times 1634509733 \ x^{12}  ) } $$

$$ \eqalign {
& ch(\theta_{21}) = { 7 \over 2^{24} \times 3^{12} \times 5^5 \ x^6 } \ (  \cr
& 2^3 \times 11 \times 1091 \times 12756480773 +  \cr
& 3 \times 5 \times 7 \times 4954321 \times 8800968989 \ x \ +  \cr
& 3^2 \times 7 \times 106275811 \times 430270161227 \ x^2 \ +  \cr
& 2^3 \times 3 \times 5 \times 7^2 \times 13 \times 3137 \times
25147 \times 246251 \times 502781 \times \ x^3 \ +  \cr
& 2^2 \times 3^3 \times 82567 \times 255832037 \times 34487228383 \ x^4 \ +  \cr
& 2^5 \times 3^3 \times 7 \times 23 \times 28753 \times 45439 \times
24160314328237 \ x^5 \ +  \cr
& 2^7 \times 3^5 \times 5 \times 7 \times 23 \times 31 \times 1867 \times
102273597906161063 \ x^6 \ +  \cr
& 2^9 \times 3^5 \times 5 \times 7 \times 23 \times 28753 \times 45439 \times
24160314328237 \ x^7  +  \cr
& 2^{10} \times 3^7 \times 5^2 \times 82567 \times 255832037 \times
34487228383 \ x^8 \ +  \cr
& 2^{15} \times 3^7 \times 5^4 \times 7^2 \times 13 \times 3137 \times
25147 \times 246251 \times 502781 \ x^9 \ +  \cr
& 2^{16} \times 3^{10} \times 5^4 \times 7 \times 106275811 \times
430270161227 \ x^{10} \ +  \cr
& 2^{20} \times 3^{11} \times 5^6 \times 7 \times 4954321 \times 8800968989 \ x^{11} \ +  \cr
& 2^{27} \times 3^{12} \times 5^6 \times 11 \times 1091 \times 12756480773 \ x^{12}  ) \cr
& \ \ \ \ \cr
& ch(\theta_{22}) = { 1 \over 2^{24} \times 3^{12} \times 5^6 \ x^7 } \ (  \cr
& 2^2 \times 3^2 \times 5^2 \times 7 \times 23 \times 443 \times 110339 \ +  \cr
& 3^2 \times 5 \times 7^2 \times 23 \times 443 \times 45307 \times 110339 \ x \ +  \cr
& 2 \times 5 \times 7 \times 79 \times 107166946015042819 \ x^2 \ +  \cr
& 2 \times 3 \times 7^2 \times 36713 \times 34478412780970831 \ x^3 \ +  \cr
& 3^3 \times 7 \times 8069 \times 48408629149359661453 \ x^4 \ +  \cr
& 2^3 \times 3^2 \times 13 \times 17 \times 4099 \times
26525879 \times 4080242157239 \ x^5 \ +  \cr
& 2^4 \times 3^3 \times 7 \times 127701402659006223221510827 \ x^6 \ +  \cr
& 2^6 \times 3^4 \times 5 \times 401 \times 1257177934278466121552321 \ x^7 \ +  \cr
& 2^8 \times 3^5 \times 5 \times 7 \times 127701402659006223221510827 \ x^8 \ +  \cr
& 2^{11} \times 3^6 \times 5^2 \times 13 \times 17 \times 4099 \times
26525879 \times 4080242157239 \ x^9 \ +  \cr
& 2^{12} \times 3^9 \times 5^3 \times 7 \times 8069 \times
48408629149359661453 \ x^{10} +  \cr
& 2^{17} \times 3^9 \times 5^4 \times 7^2 \times 36713 \times
34478412780970831 \ x^{11} \ +  \cr
& 2^{21} \times 3^{10} \times 5^6 \times 7 \times 79 \times
107166946015042819 \ x^{12} \ +  \cr
& 2^{24} \times 3^{14} \times 5^7 \times 7^2 \times 23 \times 443 \times
45307 \times 110339 \ x^{13} \ +  \cr
& 2^{30} \times 3^{16} \times 5^9 \times 7 \times 23 \times 443 \times
110339 \ x^{14}  ) } $$

$$ \eqalign {
& ch(\theta_{23}) = { 1 \over 2^{12} \times 3^6 \times 5^2 \ x^3 } \ (  \cr
& 2^6 \times 3 \times 5 \times 7 \times 51977 \ +  \cr
& 3^2 \times 37 \times 967 \times 991 \times 1277 \ x \ +  \cr
& 3 \times 5 \times 7^3 \times 37 \times 67 \times 199 \times 58321  x^2 \ +  \cr
& 2^3 \times 11 \times 104090304018661 \ x^3 \ +  \cr
& 2^4 \times 3^3 \times 5^2 \times 7^3 \times 37 \times 67 \times 199 \times
58321 \ x^4 \ +  \cr
& 2^8 \times 3^6 \times 5^2 \times 37 \times 967 \times 991 \times 1277 \ x^5 \ +  \cr
& 2^{18} \times 3^7 \times 5^4 \times 7 \times 51977 \ x^6) } $$

$$ \eqalign {
& ch(\theta_{24}) = { 7 \over 2^{10} \times 3^4 \times 5^2 \ x^3 } \ (  \cr
& 2^2 \times 3 \times 5 \times 17 \ +  \cr
& 3 \times 7 \times 17 \times 45307 \ x \ +  \cr
& 5 \times 7103 \times 192949 \ x^2 \ +  \cr
& 2^2 \times 5 \times 29 \times 541 \times 1456541 \ x^3 \ +  \cr
& 2^4 \times 3^2 \times 5^2 \times 7103 \times 192949 \ x^4 \ +  \cr
& 2^8 \times 3^5 \times 5^2 \times 7 \times 17 \times 45307 \ x^5 \ +  \cr
& 2^{14} \times 3^7 \times 5^4 \times 17 \ x^6) } $$

$$ \eqalign {
& ch(\theta_{25}) = { 1 \over 2^{16} \times 3^8 \times 5^3 \ x^4 } \ (  \cr
& 2 \times 7 \times 11 \times 1634509733 \ +  \cr
& 3 \times 5 \times 7 \times 11149 \times 240779587 \ x \ +  \cr
& 2^3 \times 97 \times 5019793 \times 33822727 \ x^2 \ +  \cr
& 3^3 \times 7^2 \times 19 \times 1429 \times 862136910847 \ x^3 \ +  \cr
& 2^2 \times 3^5 \times 7 \times 11 \times 127 \times 190523 \times
803051407 \ x^4 \ +  \cr
& 2^4 \times 3^5 \times 5 \times 7^2 \times 19 \times 1429 \times
862136910847 \ x^5 \ +  \cr
& 2^{11} \times 3^4 \times 5^2 \times 97 \times 5019793 \times 33822727 \ x^6 \ +  \cr
& 2^{12} \times 3^7 \times 5^4 \times 7 \times 11149 \times 240779587 \ x^7 \ +  \cr
& 2^{17} \times 3^8 \times 5^4 \times 7 \times 11 \times 1634509733 \ x^8) } $$

$$ \eqalign {
& ch(\theta_{26}) = { 1 \over 2^{20} \times 3^{10} \times 5^4 \ x^5 } \ (  \cr
& 3^2 \times 7 \times 1920917020339 \ +  \cr
& 2^6 \times 1932627887748647 \ x \ +  \cr
& 2^2 \times 3 \times 7 \times 687598205349751411 \ x^2 \ +  \cr
& 5 \times 7^2 \times 673 \times 77746343167071049 \ x^3 \ +  \cr
& 2^2 \times 3^2 \times 7 \times 41 \times 4513 \times 8209 \times
397543  \times 6033427 \ x^4 \ +  \cr
& 2^4 \times 3^2 \times 42256259 \times 5662271823471397 \ x^5 \ +  \cr
& 2^6 \times 3^4 \times 5 \times 7 \times 41 \times 4513 \times
8209 \times 397543 \times 6033427 \ x^6 \ +  \cr
& 2^8 \times 3^4 \times 5^3 \times 7^2 \times 673 \times 77746343167071049 \ x^7 \ +  \cr
& 2^{14} \times 3^7 \times 5^3 \times 7 \times 687598205349751411 \ x^8 \ +  \cr
& 2^{22} \times 3^8 \times 5^4 \times 1932627887748647 \ x^9 \ +  \cr
& 2^{20} \times 3^{12} \times 5^5 \times 7 \times 1920917020339 \ x^{10}) \cr
& ch(\theta_{27}) = {1 \over 2^{20} \times 3^9 \times 5^4 \ x^5 } \ (  \cr
& 7^2 \times 63113 \times 7700299 \ +  \cr
& 5^2 \times 7 \times 1699 \times 195399775663 \ x \ +  \cr
& 2^3 \times 7 \times 875675732223561979 \ x^2 \ +  \cr
& 2^3 \times 3 \times 103 \times 367 \times 5791 \times 125219 \times
13800671 \ x^3 \ +  \cr
& 2^4 \times 7 \times 233 \times 95531 \times 6117073 \times 43889063 \ x^4 \ +  \cr
& 2^4 \times 3^2 \times 7 \times 31 \times 790500401078791217977 \ x^5 \ +  \cr
& 2^8 \times 3^2 \times 5 \times 7 \times 233 \times 95531 \times
6117073 \times 43889063 \ x^6 \ +  \cr
& 2^{11} \times 3^5 \times 5^2 \times 103 \times 367 \times 5791 \times
125219 \times 13800671 \ x^7 \ +  \cr
& 2^{15} \times 3^6 \times 5^3 \times 7 \times 875675732223561979 \ x^8 \ +  \cr
& 2^{16} \times 3^8 \times 5^6 \times 7 \times 1699 \times 195399775663 \ x^9 \ +  \cr
& 2^{20} \times 3^{10} \times 5^5 \times 7^2 \times 63113 \times 7700299 x^{10} ) \cr
& ch(\theta_{28}) = { 7 \over 2^{12} \times 3^6 \times 5^2 \ x^3 } \ (  \cr
& 3 \times 7 \times 41 \times 13469 \ +  \cr
& 5 \times 25841 \times 172357 \ x \ +  \cr
& 2^2 \times 3 \times 7 \times 17 \times 151 \times 229 \times 76837 \ x^2 \ +  \cr
& 2^2 \times 3^2 \times 13 \times 41 \times 8271569177 \ x^3 \ +  \cr
& 2^6 \times 3^3 \times 5 \times 7 \times 17 \times 151 \times 229 \times
76837 \ x^4 \ +  \cr
& 2^8 \times 3^4 \times 5^3 \times 25841 \times 172357 \ x^5 \ +  \cr
& 2^{12} \times 3^7 \times 5^3 \times 7 \times 41 \times 13469 \ x^6 )   } $$

$$ \eqalign {
& ch(\theta_{29}) = { 1 \over 2^{22} \times 3^{12} \times 5^5 \ x^6 } \ (  \cr
& 3^2 \times 7 \times 37 \times 87365639293 \ +  \cr
& 3 \times 7 \times 17 \times 25158739 \times 163675601 \ x \ +  \cr
& 2 \times 7^2 \times 293207 \times 571321 \times 153212833 \ x^2 \ +  \cr
& 2 \times 3 \times 5^2 \times 7^2 \times 19 \times 4966016564031243349 \ x^3 \ +  \cr
& 2^2 \times 3^2 \times 79 \times 383 \times 1049 \times
29339 times 2479950110851 \ x^4 \ +  \cr
& 2^2 \times 3^3 \times 7 \times 259650710821 \times 25089495454453  x^5 \ +  \cr
& 2^7 \times 3^4 \times 7^2 \times 333504308440846597412509 \ x^6 \ +  \cr
& 2^6 \times 3^5 \times 5 \times 7 \times 259650710821 \times
25089495454453 \ x^7 \ +  \cr
& 2^{10} \times 3^6 \times 5^2 \times 79 \times 383 \times 1049 \times
29339 times 2479950110851 \ x^8 \ +  \cr
& 2^{13} \times 3^7 \times 5^5 \times 7^2 \times 19 \times
4966016564031243349 \ x^9 \ +  \cr
& 2^{17} \times 3^8 \times 5^4 \times 7^2 \times 293207 \times 571321 \times
153212833 \ x^{10} \ +  \cr
& 2^{20} \times 3^{11} \times 5^5 \times 7 \times 17 \times 25158739 \times
163675601 \ x^{11} \ +  \cr
& 2^{24} \times 3^{14} \times 5^6 \times 7 \times 37 \times
87365639293 \ x^{12} )    } $$

$$ \eqalign {
& ch(\theta_{30}) = { 1 \over 2^{20} \times 3^{10} \times 5^5 \ x^6 } \ (  \cr
& 2^4 \times 3^2 \times 5^3 \times 25841 +  \cr
& 2^2 \times 3^2 \times 5^2 \times 7 \times 25841 \times 45307 \ x \ +  \cr
& 5 \times 7 \times 25841 \times 115259 \times 666353 \ x^2 \ +  \cr
& 3 \times 7 \times 55717 \times 2804293 \times 13315567 \ x^3 \ +  \cr
& 2^3 \times 3^2 \times 86399 \times 959723737277659 \ x^4 \ +  \cr
& 2^6 \times 3^3 \times 5 \times 7 \times 19 \times 372825737 \times
955373927 \ x^5 \ +  \cr
& 2^8 \times 3^4 \times 5 \times 11 \times 59 \times 1627 \times
137055035733221 \ x^6 \ +  \cr
& 2^{10} \times 3^5 \times 5^2 \times 7 \times 19 \times 372825737 \times
955373927 \ x^7 \ +  \cr
& 2^{11} \times 3^6 \times 5^2 \times 86399 \times 959723737277659 \ x^8 \ +  \cr
& 2^{12} \times 3^7 \times 5^3 \times 7 \times 55717 \times 2804293 \times
13315567 \ x^9 \ +  \cr
& 2^{16} \times 3^8 \times 5^5 \times 7 \times 25841 \times 115259 \times
666353 \ x^{10} \ +  \cr
& 2^{22} \times 3^{12} \times 5^7 \times 7 \times 25841 \times 45307 \ x^{11} \ +  \cr
& 2^{28} \times 3^{14} \times 5^9 \times 25841 \ x^{12}  )   } $$

$$ \eqalign {
& ch(\theta_{31}) = { 1 \over 2^{20} \times 3^{10} \times 5^5 \ x^6 } \ (  \cr
& 2^2 \times 3^2 \times 5 \times 71 \times 1558283 +  \cr
& 3^2 \times 7 \times 71 \times 45307 \times 1558283 \ x \ +  \cr
& 2 \times 3 \times 7^2 \times 1499 \times 3758088575183 \ x^2 \ +  \cr
& 2 \times 5^2 \times 7^2 \times 11261 \times 185753 \times 172860353 \ x^3 \ +  \cr
& 3^2 \times 5^2 \times 13 \times 29 \times 241069 \times 7378760342027 \ x^4 \ +  \cr
& 2^3 \times 3^3 \times 5 \times 7 \times 11 \times 136666703950440834769 \ x^5 \ +  \cr
& 2^5 \times 3^4 \times 5^2 \times 7^2 \times 1216883197 \times 110293780367 \ x^6 \ +  \cr
& 2^7 \times 3^5 \times 5^2 \times 7 \times 11 \times 136666703950440834769 \ x^7 \ +  \cr
& 2^8 \times 3^6 \times 5^4 \times 13 \times 29 \times 241069 \times
7378760342027 \ x^8 \ +  \cr
& 2^{13} \times 3^6 \times 5^5 \times 7^2 \times 11261 \times 185753 \times
172860353 \ x^9 \ +  \cr
& 2^{17} \times 3^9 \times 5^4 \times 7^2 \times 1499 \times 3758088575183 \ x^{10} \ +  \cr
& 2^{20} \times 3^{12} \times 5^5 \times 7 \times 71 \times 45307 \times 1558283 \ x^{11} \ +  \cr
& 2^{26} \times 3^{14} \times 5^7 \times 71 \times 1558283 \ x^{12}  )   } $$

$$ \eqalign {
& ch(\theta_{32}) = { 7 \over 2^{16} \times 3^8 \times 5^4 \ x^5 } \ ( \cr
& 2^4 \times 3^3 \times 5^2 \times 17 \ +  \cr
& 2^2 \times 3^3 \times 5 \times 7 \times 17 \times 45307 \ x \ +  \cr
& 3 \times 7 \times 17 \times 115259 \times 666353 \ x^2 \ +  \cr
& 5 \times 2512868354279147 \ x^3 \ +  \cr
& 2^2 \times 3^2 \times 5 \times 1061 \times 218249 \times 31482709 \ x^4 \ +  \cr
& 2^4 \times 3^2 \times 5^2 \times 7 \times 2069 \times 634759 \times
1884341 \ x^5 \ +  \cr
& 2^6 \times 3^4 \times 5^2 \times 1061 \times 218249 \times 31482709 \ x^6 \ +  \cr
& 2^8 \times 3^4 \times 5^3 \times 2512868354279147 \ x^7 \ +  \cr
& 2^{12} \times 3^7 \times 5^3 \times 7 \times 17 \times 115259 \times
666353 \ x^8 \ +  \cr
& 2^{18} \times 3^{11} \times 5^5 \times 7 \times 17 \times 45307 \ x^9 \ +  \cr
& 2^{24} \times 3^{13} \times 5^7 \times 17 \ x^{10}  )   } $$
$$ \eqalign {
& ch(\theta_{33}) = { 1 \over 2^{14} \times 3^6 \times 5^3 \ x^4 } \ ( \cr
& 2^2 \times 5^2 \times 25841 \ +  \cr
& 5 \times 7 \times 25841 \times 45307 \ x \ +  \cr
& 3 \times 7^2 \times 19 \times 83 \times 2459 \times 45307 \ x^2 \ +  \cr
& 2^3 \times 3 \times 7^2 \times 43 \times 53841866849 \ x^3 \ +  \cr
& 2^5 \times 3^4 \times 5 \times 71 \times 109 \times 401 \times
1543 \times 1789 \ x^4 \ +  \cr
& 2^7 \times 3^3 \times 5 \times 7^2 \times 43 \times 53841866849  x^5 \ +  \cr
& 2^8 \times 3^5 \times 5^2 \times 7^2 \times 19 \times 83 \times 2459 \times
45307 \ x^6 \ +  \cr
& 2^{12} \times 3^6 \times 5^4 \times 7 \times 25841 \times 45307 \ x^7 \ +  \cr
& 2^{18} \times 3^8 \times 5^6 \times 25841 \ x^8  )  \cr
& ch(\theta_{34}) =  {1 \over 2^{20} \times 3^{10} \times 5^4 \ x^5 } \ ( \cr
& 2^5 \times 3 \times 7 \times 19 \times 101 \times 2437 \times 7883 +  \cr
& 3^3 \times 7 \times 17 \times 43 \times 670034705677 \ x \ +  \cr
& 3 \times 5 \times 7^2 \times 1093 \times 72068616146713 \ x^2 \ +  \cr
& 2^2 \times 5 \times 7^2 \times 144593 \times 92161130410753 \ x^3 \ +  \cr
& 2^2 \times 3^3 \times 7 \times 4127 \times 8389 \times 31543 \times
1323118663 \ x^4 \ +  \cr
& 2^7 \times 3^2 \times 11 \times 11943616469 \times 264210902899 \ x^5 \ +  \cr
& 2^6 \times 3^5 \times 5 \times 7 \times 4127 \times 8389 \times
31543 \times 1323118663 \ x^6 \ +  \cr
& 2^{10} \times 3^4 \times 5^3 \times 7^2 \times 144593 \times
92161130410753 \ x^7 \ +  \cr
& 2^{12} \times 3^7 \times 5^4 \times 7^2 \times 1093 \times
72068616146713 \ x^8 \ +  \cr
& 2^{16} \times 3^{11} \times 5^4 \times 7 \times 17 \times 43 \times
670034705677 \ x^9 \ +  \cr
& 2^{25} \times 3^{11} \times 5^5 \times 7 \times 19 \times 101 \times
2437 \times 7883 \ x^{10}  )  \cr
& ch(\theta_{35}) = { 1  \over 2^{14} \times 3^8 \times 5^3 \ x^4 } \ ( \cr
& 3^2 \times 71 \times 1558283 \ +  \cr
& 3 \times 7^2 \times 11 \times 71 \times 57706391 \ x \ +  \cr
& 5 \times 7 \times 71 \times 569 \times 15083 \times 158143 \ x^2 \ +  \cr
& 3^2 \times 5 \times 7^2 \times 11 \times 76091 \times 256751687 \ x^3 \ +  \cr
& 2^6 \times 3^2 \times 5 \times 238247 \times 29004591751 \ x^4 \ +  \cr
& 2^4 \times 3^4 \times 5^2 \times 7^2 \times 11 \times 76091 \times
256751687 \ x^5 \ +  \cr
& 2^8 \times 3^4 \times 5^3 \times 7 \times 71 \times 569 \times 15083 \times
158143 \ x^6 \ +  \cr
& 2^{12} \times 3^7 \times 5^3 \times 7^2 \times 11 \times 71 \times 57706391 \ x^7 \ +  \cr
& 2^{16} \times 3^{10} \times 5^4 \times 71 \times 1558283 \ x^8  )  } $$
$$ \eqalign {
& ch(\theta_{36}) = { 1  \over 2^{18} \times 3^8 \times 5^4 \ x^5} \ ( \cr
& 2^8 \times 3 \times 5^2 \times 7 \times 51977 \ +  \cr
& 2^6 \times 3 \times 5 \times 7^2 \times 45307 \times 51977 \  x \ +  \cr
& 3 \times 7 \times 13 \times 577 \times 696629 \times 1185337 \  x^2 \ +  \cr
& 5 \times 509 \times 547 \times 37532952421247 \  x^3 \ +  \cr
& 2^3 \times 7 \times 255247 \times 401008447405121 \  x^4 \ +  \cr
& 2^4 \times 5 \times 53 \times 2861 \times 1461797 \times 14111392523 \  x^5 \ +  \cr
& 2^7 \times 3^2 \times 5 \times 7 \times 255247 \times 401008447405121 \  x^6 \ +  \cr
& 2^8 \times 3^4 \times 5^3 \times 509 \times 547 \times 37532952421247 \  x^7  +  \cr
& 2^{12} \times 3^7 \times 5^3 \times 7 \times 13 \times 577 \times
696629 \times 1185337 \  x^8 \ +  \cr
& 2^{22} \times 3^9 \times 5^5 \times 7^2 \times 45307 \times 51977 \  x^9 \ +  \cr
& {2^{28} \times 3^{11} \times 5^7 \times 7 \times 51977^1} \  x^{10}  )  } $$

$$ \eqalign {
& ch(\theta_{37}) = { 1   \over 2^{22} \times 3^{10} \times 5^5 \ x^6 } \ ( \cr
& 2^2 \times 5 \times 7 \times 31 \times 193 \times 3624787 \ +  \cr
& 7^2 \times 31 \times 193 \times 45307 \times 3624787 \ x \ +  \cr
& 2 \times 31 \times 79 \times 6043 \times 672439 \times 3360061 \ x^2 \ +  \cr
& 3 \times 5 \times 7^2 \times 40277 \times 1031399 \times 973746491 \  x^3 \ +  \cr
& 2^2 \times 3 \times 5 \times 7 \times 29 \times 43 \times 97 \times 1375243 \times
57552822853 \ x^4 \ +  \cr
& 2^4 \times 3 \times 7^2 \times 239 \times 45007 \times 10354439038651493 \  x^5 \ +  \cr
& 2^6 \times 3^6 \times 5 \times 11 \times 107 \times 163 \times
208977919833364987 \ x^6 \ +  \cr
& 2^8 \times 3^3 \times 5 \times 7^2 \times 239 \times 45007 \times
10354439038651493 \  x^7 \ +  \cr
& 2^{10} \times 3^5 \times 5^3 \times 7 \times 29 \times 43 \times 97 \times
1375243 \times 57552822853 \  x^8 \ +  \cr
& 2^{12} \times 3^7 \times 5^4 \times 7^2 \times 40277 \times 1031399 \times
973746491 \  x^9 \ +  \cr
& 2^{17} \times 3^8 \times 5^4 \times 31 \times 79 \times 6043 \times
672439 \times 3360061 \ x^{10} \ +  \cr
& 2^{20} \times 3^{10} \times 5^5 \times 7^2 \times 31 \times 193 \times
45307 \times 3624787 \  x^{11} \ +  \cr
& 2^{26} \times 3^{12} \times 5^7 \times 7 \times 31 \times 193 \times
3624787 \ x^{12}  ) } $$
$$ \eqalign {
& ch(\theta_{38}) = { 7 \over 2^{16} \times 3^8 \times 5^3 \ x^4 } \ ( \cr
& 31 \times 193 \times 3624787 \ +  \cr
& 2 \times 3 \times 661 \times 39827 \times 198031 \  x \ +  \cr
& 3^4 \times 5 \times 9371 \times 3530299103 \  x^2 \ +  \cr
& 2^2 \times 3 \times 5 \times 91297 \times 232912445897 \  x^3 \ +  \cr
& 2^4 \times 3^3 \times 370597 \times 315220711859 \  x^4 \ +  \cr
& 2^6 \times 3^3 \times 5^2 \times 91297 \times 232912445897 \  x^5 \ +  \cr
& 2^8 \times 3^8 \times 5^3 \times 9371 \times 3530299103 \  x^6 \ +  \cr
& 2^{13} \times 3^7 \times 5^3 \times 661 \times 39827 \times 198031 \  x^7 \ +  \cr
& 2^{16} \times 3^8 \times 5^4 \times 31 \times 193 \times 3624787 \  x^8  )  } $$

$$ \eqalign {
& ch(\theta_{39}) = { 1 \over 2^{18} \times 3^8 \times 5^4 \times x^5 } \ ( \cr
& 2^2 \times 3 \times 5 \times 7^2 \times 41 \times 13469 \ +  \cr
& 3 \times 7^3 \times 41 \times 13469 \times 45307 \ x \ +  \cr
& 5 \times 7 \times 19813 \times 26321 \times 2574193 \ x^2 \ +  \cr
& 2^2 \times 3 \times 7 \times 41946929 \times 2853739787 \ x^3 \ +  \cr
& 2^2 \times 3 \times 7 \times 11 \times 875056131284689657 \ x^4 \ +  \cr
& 2^4 \times 3^2 \times 5^2 \times 8352802228682522537 \ x^5 \ +  \cr
& 2^6 \times 3^3 \times 5 \times 7 \times 11 \times 875056131284689657 \ x^6 \ +  \cr
& 2^{10} \times 3^5 \times 5^2 \times 7 \times 41946929 \times 2853739787 \ x^7 \ +  \cr
& 2^{12} \times 3^6 \times 5^4 \times 7 \times 19813 \times 26321 \times
2574193 \ x^8 \ +  \cr
& 2^{16} \times 3^9 \times 5^4 \times 7^3 \times 41 \times 13469 \times 45307 \ x^9 \ +  \cr
& 2^{22} \times 3^{11} \times 5^6 \times 7^2 \times 41 \times 13469 \ x^{10}  )  } $$

$$ \eqalign {
& ch(\theta_{40}) = { 7  \over 2^{18} \times 3^{10} \times 5^4 \times x^5 } \ ( \cr
& 3^2 \times 5 \times 23 \times 443 \times 110339 \ +  \cr
& 5 \times 29 \times 31 \times 113 \times 25841 \times 25849 \ x \ +  \cr
& 3 \times 7 \times 19 \times 29 \times 31 \times 83 \times 113 \times
2459 \times 25849 \ x^2 \  +  \cr
& 3^2 \times 7 \times 419 \times 2364953 \times 600815491 \ x^3 \ +  \cr
& 2^{10} \times 3^3 \times 11 \times 23 \times 8912381 \times 41333989 \ x^4 \ +  \cr
& 2^4 \times 3^4 \times 61 \times 67 \times 17529680553917933 \ x^5 \ +  \cr
& 2^{14} \times 3^5 \times 5 \times 11 \times 23 \times 8912381 \times
41333989 \ x^6 \ +  \cr
& 2^8 \times 3^6 \times 5^2 \times 7 \times 419 \times 2364953 \times 600815491 \ x^7 \ +  \cr
& 2^{12} \times 3^7 \times 5^3 \times 7 \times 19 \times 29 \times
31 \times 83 \times 113 \times 2459 \times 25849 \ x^8 \ +  \cr
& 2^{16} \times 3^8 \times 5^5 \times 29 \times 31 \times 113 \times 25841 \times
25849 \ x^9 \ +  \cr
& 2^{20} \times 3^{12} \times 5^6 \times 23 \times 443 \times 110339 \ x^{10}  )  } $$
Now, one can see that the characters above fulfill the following equation:
$$ \eqalign {
ch(\Lambda_3) \times ch(\Lambda_4) = &ch(\Lambda_3 + \Lambda_4) \ + \cr
& ch(\theta_2) \ + \ ch(\theta_3)	\ + \ ch(\theta_4) \ + \ ch(\theta_5) \ + \cr
& ch(\theta_6) \ + \ ch(\theta_7) \ + \ ch(\theta_8) \ + \ ch(\theta_9)	\ + \cr
& ch(\theta_{10}) \ + \ ch(\theta_{11})	\ + \ ch(\theta_{12}) \ + \ ch(\theta_{13})	\ + \cr
& ch(\theta_{14}) \ + \ ch(\theta_{15}) \ + \ ch(\theta_{16}) \ + \ ch(\theta_{17}) \ + \cr
& 2 \ ch(\theta_{18}) \ + \ 2 \ ch(\theta_{19}) \ + \ 2 \ ch(\theta_{20}) \ + \
2 \ ch(\theta_{21}) \ + \ 2 \ ch(\theta_{22}) \ + \cr
& 3 \ ch(\theta_{23}) \ + \ 3 \ ch(\theta_{24}) \ + \ 3 \ ch(\theta_{25}) \ +
3 \ ch(\theta_{26}) \ +  \cr
& 3 \ ch(\theta_{27}) \ + \ 3 \ ch(\theta_{28}) \ + \ 3 \ ch(\theta_{29}) \ + \ 3 \ ch(\theta_{30}) \ + \cr
& 4 \ ch(\theta_{31}) \ + \ 4 \ ch(\theta_{32}) \ + \ 4 \ ch(\theta_{33}) \ +  \cr
& 5 \ ch(\theta_{34}) \ + \ 5 \ ch(\theta_{35}) \ + \ 5 \ ch(\theta_{36}) \ + \ 5 \ ch(\theta_{37}) \ +  \cr
& 5 \ ch(\theta_{38}) \ + \ 5 \ ch(\theta_{39}) \ + \ 5 \ ch(\theta_{40})	}   \eqno(III.6)  $$

\noindent One should note however that, the 1-parameter specialization (III.5) above is not enough to find all the tensor coupling coefficients completely so we saw that at least 3-parameters specializations will be sufficient, which we used the following one;
$$ \eqalign {
u_1 &= 1 \cr
u_2 &= 2 \cr
u_3 &= 3 \cr
u_4 &= 4 \cr
u_5 &= x \cr
u_6 &= y \cr
u_7 &= z \cr
u_8 &= 1/(24 \ x \ y \ z)  }  $$
\vskip 3mm
\noindent {\bf{IV.\ APPENDIX }}
$$ \eqalign {
&\Sigma(1)=1 \cr &\Sigma(2)=\sigma_{3,2} \ , \
\Sigma(3)=\sigma_{3,4} \cr &\Sigma( 4)=\sigma_{3,2,1,4} \  , \
\Sigma( 5)=\sigma_{3,2,4,3} \ , \ \Sigma( 6)=\sigma_{3,2,4,5} \cr
&\Sigma( 7)=\sigma_{3,2,4,7} \ , \ \Sigma( 8)=\sigma_{3,4,5,6} \ , \
\Sigma( 9)=\sigma_{3,4,5,7} \cr &\Sigma(10)=\sigma_{3,2,1,4,3,2} \ ,
\ \Sigma(11)=\sigma_{3,2,1,4,3,5} \ , \
\Sigma(12)=\sigma_{3,2,1,4,3,7} \cr &\Sigma(13)=\sigma_{3,2,1,4,5,6}
\, \ \Sigma(14)=\sigma_{3,2,1,4,5,7} \ , \
\Sigma(15)=\sigma_{3,2,4,3,5,4} \cr &\Sigma(16)=\sigma_{3,2,4,3,5,6}
\ , \ \Sigma(17)=\sigma_{3,2,4,3,5,7} \ , \
\Sigma(18)=\sigma_{3,2,4,3,7,4} \cr &\Sigma(19)=\sigma_{3,2,4,5,6,7}
\ , \ \Sigma(20)=\sigma_{3,2,4,5,7,4} \ , \
\Sigma(21)=\sigma_{3,4,5,6,7,4} \cr &\Sigma(22)=\sigma_{3,4,5,7,4,3}
\cr &\Sigma(23)=\sigma_{3,2,1,4,3,2,5,4} \ , \
\Sigma(24)=\sigma_{3,2,1,4,3,2,5,6} \ , \
\Sigma(25)=\sigma_{3,2,1,4,3,2,5,7} \cr
&\Sigma(26)=\sigma_{3,2,1,4,3,2,7,4} \ , \
\Sigma(27)=\sigma_{3,2,1,4,3,5,4,6}  \ , \
\Sigma(28)=\sigma_{3,2,1,4,3,5,6,7} \cr
&\Sigma(29)=\sigma_{3,2,1,4,5,6,7,4}  \ , \
\Sigma(30)=\sigma_{3,2,4,3,5,4,6,5}  \ , \
\Sigma(31)=\sigma_{3,2,4,5,6,7,4,5} \cr
&\Sigma(32)=\sigma_{3,4,5,6,7,4,3,5} \cr
&\Sigma(33)=\sigma_{3,2,1,4,3,2,5,4,3,6}  \ , \
\Sigma(34)=\sigma_{3,2,1,4,3,2,5,4,6,5}  \ , \
\Sigma(35)=\sigma_{3,4,5,6,7,4,3,5,4,7}  \cr
&\Sigma(36)=\sigma_{3,2,1,4,3,2,5,4,3,6,5,4} \cr} $$

$$ \eqalign {
&\Sigma(37)=\sigma_3  \cr &\Sigma(38)=\sigma_{3,2,1}  \ , \
\Sigma(39)=\sigma_{3,2,4}  \ , \ \Sigma(40)=\sigma_{3,4,5}  \ , \
\Sigma(41)=\sigma_{3,4,7}  \cr &\Sigma(42)=\sigma_{3,2,1,4,3}  \ , \
\Sigma(43)=\sigma_{3,2,1,4,5} \ , \ \Sigma(44)=\sigma_{3,2,1,4,7}
\cr &\Sigma(45)=\sigma_{3,2,4,3,5}  \ , \
\Sigma(46)=\sigma_{3,2,4,3,7} \ , \ \Sigma(47)=\sigma_{3,2,4,5,6}
\cr &\Sigma(48)=\sigma_{3,2,4,5,7} \ , \
\Sigma(49)=\sigma_{3,4,5,6,7} \ , \ \Sigma(50)=\sigma_{3,4,5,7,4} \cr
&\Sigma(51)=\sigma_{3,2,1,4,3,2,5} \ , \
\Sigma(52)=\sigma_{3,2,1,4,3,2,7}  \ , \
\Sigma(53)=\sigma_{3,2,1,4,3,5,4}  \cr
&\Sigma(54)=\sigma_{3,2,1,4,3,5,6}  \ , \
\Sigma(55)=\sigma_{3,2,1,4,3,5,7}  \ , \
\Sigma(56)=\sigma_{3,2,1,4,3,7,4} \ , \cr
&\Sigma(57)=\sigma_{3,2,1,4,5,6,7}  \ , \
\Sigma(58)=\sigma_{3,2,1,4,5,7,4}  \ , \
\Sigma(59)=\sigma_{3,2,4,3,5,4,6}  \cr
&\Sigma(60)=\sigma_{3,2,4,3,5,6,7}  \ , \
\Sigma(61)=\sigma_{3,2,4,5,6,7,4}  \ , \
\Sigma(62)=\sigma_{3,4,5,6,7,4,3}  \cr
&\Sigma(63)=\sigma_{3,4,5,6,7,4,5}  \cr
&\Sigma(64)=\sigma_{3,2,1,4,3,2,5,4,3}  \ , \
\Sigma(65)=\sigma_{3,2,1,4,3,2,5,4,6}  \ , \
\Sigma(66)=\sigma_{3,2,1,4,3,2,5,6,7}  \cr
&\Sigma(67)=\sigma_{3,2,1,4,3,2,7,4,3}  \ , \
\Sigma(68)=\sigma_{3,2,1,4,3,5,4,6,5}  \ , \
\Sigma(69)=\sigma_{3,2,1,4,5,6,7,4,5}  \cr
&\Sigma(70)=\sigma_{3,4,5,6,7,4,3,5,4}  \cr
&\Sigma(71)=\sigma_{3,2,1,4,3,2,5,4,3,6,5}  \cr
&\Sigma(72)=\sigma_{3,2,1,4,3,2,5,4,3,6,5,4,7}  }  \eqno(II.6) $$

\vskip 3mm
\noindent {\bf{V.\ REFERENCES }}
\vskip 3mm

\item [1] J. E. Humphreys, Introduction to Lie Algebras and Representation Theory,
\item \ \ \ \ Springer-Verlag, 1972
\item [2] 24.4 in ref.[1]
\item [3] V. Kac, Infinite Dimensional Lie Algebras, Cambridge University Press, 1982
\item [4] H.R.Karadayi, M.Gungormez, Fundamental Weights, Permutation Weights and
\item \ \ \ \ \ Weyl Character Formula, J.Phys.A32:1701-1707(1999)

\end